\documentclass[11pt]{article}
\usepackage{amsfonts}
\usepackage{amssymb, amsmath, cite}

\newtheorem{theorem}{Theorem}
\newtheorem{lemma}{Lemma}
\newtheorem{proposition}{Proposition}
\newtheorem{remark}{Remark}

\newcommand\be{\begin{equation}}
\newcommand\ee{\end{equation}}
\newcommand\ber{\begin{eqnarray}}
\newcommand\eer{\end{eqnarray}}
\newcommand\berr{\begin{eqnarray*}}
\newcommand\eerr{\end{eqnarray*}}
\newcommand{\ud}{\mathrm{d}}\newcommand{\bfR}{\mathbb{R}}\newcommand{\pa}{\partial}\newcommand{\ii}{\mbox{i}}
\newcommand{\nm}{\nonumber}\newcommand{\bfZ}{\mathbb{Z}}\newcommand{\e}{\mbox{e}}\newcommand{\lm}{\lambda}
\newcommand{\ito}{\int_{\Omega}}
\newcommand{\itr}{\int_{\mathbb{R}^2}}
\newcommand{\vep}{\varepsilon}
 \newcommand\re{\mathrm{e}}
 \newcommand\wot  {W^{1, 2}(\mathbb{R}^2)}
{\setlength\arraycolsep{2pt}

\begin{document}

\title{Chern--Simons Vortices in the Gudnason Model}
\author{Xiaosen Han\\Institute of Contemporary Mathematics\\School of Mathematics\\Henan University\\
Kaifeng, Henan 475004, PR China\\ \\
Chang-Shou Lin\\Department of Mathematics\\ National Taiwan University\\ Taipei, Taiwan 10617, ROC\\ \\Gabriella Tarantello\\Dipartimento di Matematica\\Unversit\`{a} di
Roma ``Tor Vergata"\\Via della Ricerca Scientifica\\00133 Rome, Italy\\ \\Yisong Yang\\Department of Mathematics\\Polytechnic Institute of New York University\\Brooklyn, New York 11201, USA}

\date{}

\maketitle

\begin{abstract}
We present a series of existence theorems for multiple vortex solutions in the Gudnason model of the ${\cal N}=2$ supersymmetric field theory where non-Abelian
gauge fields are governed by the pure Chern--Simons dynamics at dual levels and realized as the solutions
of a system of elliptic equations with exponential nonlinearity over two-dimensional domains. In the full plane situation, our method utilizes a minimization approach, and in the doubly periodic situation,
we employ an-inequality constrained minimization approach. In the latter case, we also obtain sufficient conditions under which we show that there exist at least two gauge-distinct solutions for
any  prescribed distribution of vortices. In other words, there are distinct solutions with identical vortex distribution, energy, and electric and magnetic charges.
\end{abstract}

\section{Introduction}

The classical Abelian Higgs model defined by the Lagrangian action density
\be \label{x1}
{\cal L}_{\mbox{AH}}=-\frac14 F_{\mu}F^{\mu\nu}+\frac12D_\mu\phi\overline{D^\mu\phi}-\frac\lm8(|\phi|^2-1)^2,
\ee
is of foundational importance in theoretical physics. Here $\phi$ is a complex-valued scalar field called the Higgs field, $F_{\mu\nu}=\partial_\mu A_\nu-\partial_\nu A_\mu$ the
electromagnetic field strength generated from a real-valued gauge field, $A_\mu$, $D_\mu\phi=\partial_\mu\phi-\ii A_\mu\phi$ the gauge-covariant derivative, $\lm>0$ a coupling
parameter, $\mu,\nu=0,1,2,3$ are the $(3+1)$-dimensional Minskowski spacetime coordinate indices, and the metrics $(g^{\mu\nu})=(g_{\mu\nu})=\mbox{diag}(1,-1,-1,-1)$ are
used to raise or lower indices. For example, in quantum field theory, the model provides a mathematically simplest thought-laboratory allowing various
fundamental concepts such as spontaneous symmetry-breaking, the onset and annihilation of the Goldstone particles, and the Higgs mechanism to
be formulated and explored \cite{Chaichian,Huang,Ryder}. Phenomenologically, the model in its static limit and temporal gauge $A_0=0$ is the celebrated Ginzburg--Landau theory \cite{GL} for superconductivity. In two spatial dimensions, topological defects
in the form of the Abrikosov vortices \cite{Ab,Du} can be generated from the model which, when coupled with the Einstein equations,  provide indispensable structures,  known as cosmic strings, giving
rise to centers of curvature concentrations in spacetime and hence forming the seeds for matter accretion in the early universe \cite{Vi,Yob1,Yob2}. On the other hand, however, a statement known as the Julia--Zee theorem
\cite{JZ,SYjz} says
that the static solutions of the equations of motion of (\ref{x1}) as well as of that of the general non-Abelian Yang--Mills--Higgs model, in two-spatial dimensions and of finite energy, must stay in the temporal gauge, $A_0=0$. In other words, vortices of the
Yang--Mills--Higgs model, Abelian or non-Abelian, do not carry electric charge and can only be purely magnetic.  Nevertheless, electrically {\em and} magnetically charged vortices, called
dyonic vortices, are needed in many areas of applications such as high-temperature superconductivity \cite{Kh,Ma}, the Bose--Einstein condensates \cite{In,Ka}, the quantum Hall effect \cite{So},
optics \cite{Be}, and superfluids \cite{She}. Therefore, it will be
important to modify the classical Yang--Higgs--Higgs theory so that dyonic vortices are accommodated. In a series of pioneering studies,
Jackiw and Templeton \cite{JTemp}, Schonfeld \cite{Schonfeld}, Deser, Jackiw, and Templeton \cite{DJT1,DJT2}, Paul and Khare \cite{PK}, de Vega and Schaposnik \cite{VS,VS1}, and
Kumar and Khare \cite{KK} developed a modified Yang--Mills--Higgs theory in which the Chern--Simons topological terms \cite{CS1,CS2} are implemented into the action density.
Although these terms fail to be gauge-invariant locally, they preserve gauge-invariance globally and thus render the theory the same gauge invariance as in the Yang--Mills--Higgs theory. More importantly, the presence of the Chern--Simons terms make the coexistence of electric and magnetic charges a necessity. Mathematically, however, the presence of the Chern--Simons terms and the nontrivial
temporal component of the gauge field leads us to face a much more complicated form of the equations of motion and an existence theory for radially symmetric solutions
has only been obtained rather recently  \cite{CGSY}. This difficulty motivated people to explore possible BPS (after the seminal studies of Bogomol'nyi \cite{Bo} and Prasad and Sommerfield \cite{PS})
reductions of the problem and brought into light the works of Hong, Kim, and Pac \cite{HKP} and Jackiw and Weinberg \cite{JW}, which sparked a great development of the subject of construction of
the Chern--Simons--Higgs vortex solutions up to today. In their approach, it may be understood that the initial Lagrangian action density to be modified is simply that of the Abelian Higgs model, with the addition of a Chern--Simons
term controlled by a coupling parameter $\kappa$, which may be written as
\be \label{x2}
{\cal L}_{\mbox{MCSH}}=-\frac1{4e^2} F_{\mu\nu}F^{\mu\nu}-\frac\kappa4\epsilon^{\mu\nu\rho}A_\mu F_{\nu\rho}+D_\mu\phi\overline{D^\mu\phi}-V(|\phi|^2),
\ee
where $e>0$ denotes a coupling parameter imposed on the Maxwell kinetic density term and $V$ is the Higgs potential density function, as in \cite{PK,VS}.
Since in (\ref{x2}) both the Maxwell and Chern--Simons terms are present, the model is referred to as the Maxwell--Chern--Simons--Higgs model. In the limit $e\to\infty$, the
Maxwell term is switched off and the model (\ref{x2}) becomes
\be \label{x3}
{\cal L}_{\mbox{ACSH}}=-\frac\kappa4\epsilon^{\mu\nu\rho}A_\mu F_{\nu\rho}+D_\mu\phi\overline{D^\mu\phi}-V(|\phi|^2),
\ee
which is known as the Abelian Chern--Simons--Higgs model \cite{HKP,JW} in which the gauge field dynamics is governed solely by the Chern--Simons term. It is shown in \cite{HKP,JW} that, when the Higgs
potential function $V$ is chosen to be
\be
V(|\phi|^2)=\frac1{\kappa^2}|\phi|^2(1-|\phi|^2)^2,
\ee
which is analogous with the critical choice $\lm=1$ for (\ref{x1}) studied in \cite{jata},  the equations of motion of (\ref{x3}) may be reduced into the  BPS system
\be
D_1\phi+\ii D_2\phi=0,\quad F_{12}=\frac2{\kappa^2}|\phi|^2(1-|\phi|^2).
\ee
The multiple vortex solutions, realizing a prescribed distribution of vortices located at $p_1,\dots,p_n$ and carrying the total electric and magnetic charges, $Q_{\mbox{e}}$ and $Q_{\mbox{m}}$, given by
\be
Q_{\mbox{e}}=\kappa Q_{\mbox{m}},\quad Q_{\mbox{m}}=\int F_{12}\,\ud x=2\pi n,
\ee
may be constructed in terms of the variable $u=|\phi|^2$ via solving the master equation \cite{HKP,JW}
\be \label{x7}
\Delta u=\alpha\e^u(\e^u-1)+4\pi\sum_{s=1}^n\delta_{p_s}(x),
\ee
where $\alpha=\frac4{\kappa^2}$, whose structure has been shown \cite{tabook,Ybook} to be much richer and more challenging than that of the classical Abelian Higgs model. For example, in contrast to the Abelian Higgs
model (\ref{x1}) where the solution of finite energy realizing any prescribed distribution of vortices is unique \cite{jata,tau1,WY}, the solutions in the Abelian Chern--Simons--Higgs model are
not unique and further categorized \cite{Dunne1,JPW,HZ} into topological and non-topological solutions, which have led to a rapid development of analytic methods and
the harvest of a rich vista of results \cite{caya1,Chae1,Chan,CHMcY,SYcs1,SYcs2} regarding issues such as existence, uniqueness, nonexistence, asymptotic properties, approximation, etc., of the solutions.

The purpose of the present article is to develop an existence theory for non-Abelian Chern--Simon--Higgs vortices in the Gudnason model \cite{gud1,gud2}.
Historically, after the formulation of the Abelian Chern--Simons--Higgs model in \cite{HKP,JW}, a great deal of activities quickly evolved around developing  \cite{DJPT}
and analyzing \cite{LWY,nota,Ycs} non-Abelian extensions of
the model, both non-relativistic and relativistic. What distinguish the Gudnason model \cite{gud1,gud2} from the classical non-Abelian Chern--Simons--Higgs models  \cite{DJPT,VS,VS1}
are that the former may be derived in a supersymmetric gauge field theory framework as those in the study of monopole confinement mechanism \cite{Auzzi,HSZ,HT,MY,ShY2,ShY} and that
a kind of bi-level Chern--Simons
dynamics is present as in the Bagger--Lambert--Gustavsson theory \cite{BL1,BL2,BL3,G} and  the Aharony--Bergman--Jafferis--Maldacena theory \cite{abjm} to govern gauge-field kinematics.
These studies have prompted a great amount of research activities over the past few years and the nonlinear partial differential equation problems unearthed offer truly rich opportunities for analysts
in exploring new techniques and ideas.

We shall present two types of results. The first type concerns the existence of solutions over the full plane subject to the boundary behavior corresponding to the asymptotic vacuum state with
a completely broken symmetry \cite{gud1,gud2}. Our method is based on a variational reformulation of the problem and a coercive minimization approach.
The second type concerns the existence, nonexistence, and multiplicity of solutions over a doubly periodic domain. The main difficulty we encounter here is the constraint problem which makes
it hard to develop a general variational method as in the full plane situation where there is no constraint to tackle. In this situation we compromise to solve the system in the special case with
two equations. Although this case is limited, the results are rich and structures are challenging. Here we extend the inequality-constrained techniques \cite{caya1,nota,taran96} to resolve the
equality-constraint difficulty and obtain existence and multiplicity results for solutions.

A brief outline of the rest of the paper is as follows. In the next section, we first review the Gudnason model \cite{gud1,gud2} and the associated non-Abelian vortex equations, also called the master equations \cite{gud1,gud2}, to be
studied in this paper. In the following section, we reduce the master equations into a system of nonlinear equations to be analyzed and state our main existence theorems regarding
the solutions. The subsequent
sections are then devoted to the proofs of these theorems by developing and utilizing variational techniques. In the last section, we briefly summarize and comment on our results.

\section{The Gudnason model}

Consider the standard Minkowski spacetime $\bfR^{2,1}$ of signature $(+--)$ and use $\mu,\nu=0,1,2$ to denote the temporal and spatial coordinate
indices. The Gudnason model \cite{gud1} is formulated as an ${\cal N}=2$ supersymmetric Yang--Mills--Chern--Simons--Higgs theory with the general gauge group
$G=U(1)\times G'$ where $G'$ is a non-Abelian simple Lie group represented by matrices. As in \cite{gud1}, use the index $a=1,\dots,\dim(G')$ to label the non-Abelian
gauge group generators and the index $0$ the Abelian one, $U(1)$. Then the gauge potential $A_\mu$ taking values in the Lie algebra of the group $G$ may be written as
$A_\mu =A^\alpha_\mu t^\alpha$ where $A^\alpha_\mu$ are real-valued vector fields and $t^\alpha$ ($\alpha=0,a$) the generators of $G$ which are normalized to satisfy
\be
t^0=\frac1{\sqrt{2N}},\quad\mbox{Tr}(t^a t^b)=\frac12\delta^{ab},
\ee
where $N$ is the dimension of the fundamental representation space of $G'$.
The gauge field strength tensor $F_{\mu\nu}$ is then given by
\be
F_{\mu\nu}=\pa_\mu A_\nu -\pa_\nu A_\mu +\ii[A_\mu,A_\nu].
\ee
Use $\phi$ to denote the Higgs scalar field in the adjoint representation of $G$ and the $N\times N_{\small\mbox{f}}$ matrix $H$ to contain $N_{\small\mbox{f}}$-flavor
matter (quark) fields in the fundamental representation of $G$. Then their gauge-covariant derivatives are
\be
{\cal D}_\mu \phi=\partial_\mu\phi+\ii [A_\mu,\phi],\quad {\cal D}_\mu H=\partial_\mu H+\ii A_\mu H.
\ee
With this preparation, the Lagrangian action density of the Yang--Mills--Chern--Simons--Higgs theory of the Gudnason model \cite{gud1}, omitting the Fermion part, is written as
\ber\label{G1}
{\cal L}_{\small\mbox{YMCSH}}&=&-\frac1{4g^2}(F^a_{\mu\nu})^2-\frac1{4e^2}(F^0_{\mu\nu})^2-\frac{\mu}{8\pi}\epsilon^{\mu\nu\rho}\left(A^a_\mu\partial_\nu A^a_\rho-\frac13f^{abc}A^a_\mu A^b_\nu A^c_\rho
\right)-\frac\kappa{8\pi}\epsilon^{\mu\nu\rho} A^0_\mu \partial_\nu A^0_\rho\nm\\
&&+\frac1{2g^2}({\cal D}_\mu\phi^a)^2+\frac1{2e^2}(\partial_\mu\phi^0)^2+\mbox{Tr}({\cal D}_\mu H)({\cal D}^\mu H)^\dagger-\mbox{Tr}|\phi H-Hm|^2\nm\\
&&-\frac{g^2}2\left(\mbox{Tr}(H H^\dagger t^a)-\frac\mu{4\pi}\phi^a\right)^2-\frac{e^2}2\left(\mbox{Tr}(HH^\dagger t^0)-\frac\kappa{4\pi}\phi^0-\frac\xi{\sqrt{2N}}\right)^2,
\eer
where $\epsilon^{\mu\nu\nu}$ is the Kronecker skew-symmetric tensor, $f^{abc}$ are the structural constants of the non-Abelian gauge group $G'$,
$e,g\in\bfR$ the Abelian and non-Abelian Yang--Mills coupling constants, $\kappa\in\bfR$ the Abelian Chern--Simons coupling, $\mu\in\bfZ$
the non-Abelian Chern--Simons coupling constant which should not be confused with the spacetime coordinate index, $m$ is a mass matrix,  $\xi\in\bfR$ is the so-called Fayet--Iliopoulos parameter, and
the notation $|A|^2=AA^\dagger$ is observed for a complex matrix $A$.

In the strong coupling limit, $e,g\to\infty, m=0, \kappa\neq\mu$, the Gudnason model (\ref{G1}) is shown to reduce to the following Chern--Simons--Higgs model \cite{gud1,gud2}
\ber \label{G2}
{\cal L}_{\mbox{CSH}}&=& -\frac{\mu}{8\pi}\epsilon^{\mu\nu\rho}\left(A^a_\mu\partial_\nu A^a_\rho-\frac13f^{abc}A^a_\mu A^b_\nu A^c_\rho
\right)-\frac\kappa{8\pi}\epsilon^{\mu\nu\rho} A^0_\mu \partial_\nu A^0_\rho\nm\\
&&+\mbox{Tr}({\cal D}_\mu H)({\cal D}^\mu H)^\dagger-4\pi^2\mbox{Tr}\bigg|\left(\frac1{\kappa N}(\mbox{Tr}(HH^\dagger)-\xi){\bf 1}_N+\frac2\mu\mbox{Tr}(HH^\dagger t^a)t^a\right)H\bigg|^2,
\eer
where ${\bf 1}_N$ denotes the $N\times N$ identity matrix. The equations of motion of the Gudnason model \cite{gud1,gud2} governed by the Lagrangian action density (\ref{G2}) are \cite{gud1}
\ber
\frac\mu{8\pi}\epsilon^{\mu\nu\rho}F^a_{\mu\nu}&=&-\ii \mbox{Tr}\left(H^\dagger t^a {\cal D}^\rho H-({\cal D}^\rho H)^\dagger t^a H\right),\\
\frac\kappa{8\pi}\epsilon^{\mu\nu\rho}F^0_{\mu\nu}&=&-\ii \mbox{Tr}\left(H^\dagger t^0 {\cal D}^\rho H-({\cal D}^\rho H)^\dagger t^0 H\right),\\
{\cal D}_\mu {\cal D}^\mu H&=&-4\pi^2\left(\frac1{\kappa N}(\mbox{Tr}(HH^\dagger)-\xi){\bf 1}_N+\frac2\mu\mbox{Tr}(HH^\dagger t^a)t^a\right)^2 H\nm\\
&&-\frac{8\pi^2}{\kappa N}\mbox{Tr}\left\{\left(\frac1{\kappa N}(\mbox{Tr}(HH^\dagger)-\xi){\bf 1}_N+\frac2\mu\mbox{Tr}(HH^\dagger t^a)t^a\right)HH^\dagger\right\} H\nm\\
&&-\frac{16\pi^2}\mu \mbox{Tr}\left\{\left(\frac1{\kappa N}(\mbox{Tr}(HH^\dagger)-\xi){\bf 1}_N+\frac2\mu\mbox{Tr}(HH^\dagger t^b)t^b\right)HH^\dagger t^a\right\}t^a H,
\eer
which are rather complicated. To proceed, Gudnason \cite{gud1} applies the Bogomol'nyi method to carry out a completion of square analysis for the string tension of the model, obtained from the
integration of the Hamiltonian component of the stress tensor, in the static limit, and derives the following BPS system of equations:
\ber
\overline{\cal D}H&=&0,\\
F^a_{12} t^a&=&\frac{16\pi^2}{\mu\kappa N}(\mbox{Tr}(HH^\dagger)-\xi)\mbox{Tr}(HH^\dagger t^a)t^a+\frac{16\pi^2}{\mu^2}\mbox{Tr}(HH^\dagger t^b)\mbox{Tr}(HH^\dagger\{t^a,t^b\})t^a,\\
F^0_{12}t^0&=&\frac{8\pi^2}{\kappa^2 N^2}\mbox{Tr}(HH^\dagger)(\mbox{Tr}(HH^\dagger)-\xi){\bf 1}_N+\frac{16\pi^2}{\mu\kappa N}(\mbox{Tr}(HH^\dagger t^a))^2 {\bf1}_N,
\eer
where $\overline{\cal D}={\cal D}_1+\ii{\cal D}_2$. In order to investigate the solutions of these equations, Gudnason \cite{gud2} specifies the concrete situations, $G'=SO(2M)$ and $G'=USp(2M)$
with $N=2M$ so that the equations become
\ber
\overline{\cal D}H&=&0,\label{BPS1}\\
F^a_{12} t^a&=&\frac{2\pi^2}{\mu\kappa M}(\mbox{Tr}(HH^\dagger)-\xi)\langle HH^\dagger\rangle_J+\frac{2\pi^2}{\mu^2}\langle(HH^\dagger)^2\rangle_J,\label{BPS2}\\
F^0_{12}t^0&=&\frac{2\pi^2}{\kappa^2 M^2}\mbox{Tr}(HH^\dagger)(\mbox{Tr}(HH^\dagger)-\xi){\bf 1}_{2M}+\frac{2\pi^2}{\mu\kappa M}\mbox{Tr}(HH^\dagger \langle HH^\dagger\rangle_J) {\bf1}_{2M},\label{BPS3}
\eer
in which $\langle A\rangle_J=A-J^\dagger A^t J$ for a $2M\times 2M$ matrix with
\be
J=\left(\begin{array}{cc}{\bf0}&{\bf 1}_M\\\epsilon {\bf1}_M&{\bf0}\end{array}\right),\quad \epsilon=\pm1\mbox{ depending on whether }G'=SO(2M) \mbox{ or } USp(2M).
\ee
Then, choosing $H_0$ as a suitable background matrix realizing a prescribed distribution of vortices and using a moduli matrix ansatz of the form
$H=S^{-1} H_0$, $S=sS'$, which splits the variables into the Abelian one $\omega=|s|^2$ and the non-Abelian one $\Omega'=S'(S')^\dagger$, so that $\Omega=SS^\dagger=\omega\Omega'$,
the BPS equations (\ref{BPS1})--(\ref{BPS3}) are shown to become the so-called master equations \cite{gud2}
\ber
\overline{\partial}(\Omega'\partial\Omega'^{-1})&=&\frac{\pi^2}{\mu\kappa M}(\mbox{Tr}(\Omega_0\Omega^{-1})-\xi)\langle\Omega_0\Omega^{-1}\rangle_J+\frac{\pi^2}{\mu^2}\langle(\Omega_0\Omega^{-1})^2\rangle_J,\label{m1}\\
\overline{\partial}\partial\ln\omega&=&-\frac{\pi^2}{\kappa^2 M^2}\mbox{Tr}(\Omega_0\Omega^{-1})(\mbox{Tr}(\Omega_0\Omega^{-1})-\xi)-\frac{\pi^2}{\mu\kappa M}\mbox{Tr}(\Omega_0\Omega^{-1}\langle\Omega_0\Omega^{-1}\rangle_J),\label{m2}
\eer
where $\overline{\partial}=\partial_1+\ii\partial_2$.

To describe  multiple vortices  by the master equations \eqref{m1}--\eqref{m2},     we take as in \cite{gud1,gud2} the  moduli matrix  $H_0(z)$  of  the form
   \berr
 H_0(z)
   =\prod\limits_{i=1}^{M}\prod\limits_{s=1}^{n_i}(z-z_{i, s})D_0(z),
\eerr
where
\berr
D_0(z)={\rm diag}\left\{\prod\limits_{s=1}^{n_1}(z-z_{1, s}), \dots, \prod\limits_{s=1}^{n_M}(z-z_{M, s}),\prod\limits_{s=1}^{n_1}(z-z_{1, s})^{-1}, \dots, \prod\limits_{s=1}^{n_M}(z-z_{M, s})^{-1}\right\}, \label{h1}
   \eerr
 and $z_{i, s} $ are prescribed points on the complex plane, $s=1, \dots, n_i, i=1, \dots, M$; $n_i $ are nonnegative integers, $i=1, \dots, M$.
   We easily see that
      \berr
    H_0^T(z)JH_0(z)&=&\prod\limits_{i=1}^{M}\prod\limits_{s=1}^{n_i}(z-z_{i, s})J,\nm\\
    \det\left(H_0(z)\right)&=&\left(\prod\limits_{i=1}^{M}\prod\limits_{s=1}^{n_i}(z-z_{i, s})\right)^{2M},\nm\\
    \Omega_0=H_0(z)H_0(z)^\dagger
    &=&\prod\limits_{i=1}^{M}\prod\limits_{s=1}^{n_i}|z-z_{i, s}|^2 D_0(z)^\dagger D_0(z).\nm
\eerr
   With  the  further    ansatz
\be
\Omega'=\mbox{diag}(\e^{\chi_1},\dots,\e^{\chi_M},\e^{-\chi_1},\dots,\e^{-\chi_M}),\quad \omega=\e^\psi,
\ee
where $\chi_1,\dots,\chi_M,\psi$ are real-valued functions, and a direct computation,  the  equations (\ref{m1})--(\ref{m2}) read \cite{gud2}:
 \ber
\overline{\partial}\partial \chi_j&=&-\frac{\pi^2}{\mu\kappa M}\left(\prod\limits_{k=1}^{M}\prod\limits_{s=1}^{n_k}|z-z_{k, s}|^2\sum_{i=1}^M\left[\prod\limits_{s=1}^{n_i}|z-z_{i, s}|^2\e^{-\psi-\chi_i}+\prod\limits_{s=1}^{n_i}|z-z_{i, s}|^{-2}\e^{-\psi+\chi_i}\right]-\xi\right)\nm\\
&&\times\prod\limits_{k=1}^{M}\prod\limits_{s=1}^{n_k}|z-z_{k, s}|^2\left(\prod\limits_{s=1}^{n_j}|z-z_{j, s}|^2\e^{-\psi-\chi_j}-\prod\limits_{s=1}^{n_j}|z-z_{j, s}|^2\e^{-\psi+\chi_j}\right)\nm\\
&&-\frac{\pi^2}{\mu^2}\prod\limits_{k=1}^{M}\prod\limits_{s=1}^{n_k}|z-z_{k, s}|^4\left(\prod\limits_{s=1}^{n_j}|z-z_{j, s}|^4\e^{-2\psi-2\chi_j}-\prod\limits_{s=1}^{n_j}|z-z_{j, s}|^{-4}\e^{-2\psi+2\chi_j}\right),\nm\\
&&\qquad \qquad j=1,\dots,M,\label{mm1}\\
\overline{\partial}\partial\psi&=&-\frac{\pi^2}{\kappa^2 M^2}\left(\prod\limits_{k=1}^{M}\prod\limits_{s=1}^{n_k}|z-z_{k, s}|^2\sum_{i=1}^M\left[\prod\limits_{s=1}^{n_i}|z-z_{i, s}|^2\e^{-\psi-\chi_i}+\prod\limits_{s=1}^{n_i}|z-z_{i, s}|^{-2}\e^{-\psi+\chi_i}\right]-\xi\right)\nm\\
&&\times \prod\limits_{k=1}^{M}\prod\limits_{s=1}^{n_k}|z-z_{k, s}|^2\sum_{j=1}^M\left(\prod\limits_{s=1}^{n_j}|z-z_{j, s}|^2\e^{-\psi-\chi_j}+\prod\limits_{s=1}^{n_j}|z-z_{j, s}|^{-2}\e^{-\psi+\chi_j}\right)
\nm\\
&&-\frac{\pi^2}{\mu\kappa M}\sum_{i=1}^M\prod\limits_{s=1}^{n_k}|z-z_{k, s}|^4\sum_{i=1}^M\left(\prod\limits_{s=1}^{n_i}|z-z_{i, s}|^{2}\e^{-\psi-\chi_i}-\prod\limits_{s=1}^{n_i}|z-z_{i, s}|^{-2}\e^{-\psi+\chi_i}\right)^2.\label{mm2}
\eer

These are the master equations which govern the multiple vortex solutions of the  non-Abelian Chern--Simons--Higgs model of Gudnason \cite{gud1,gud2}.
 Below we aim to establish a series of existence theorems for the solutions of these equations over the full plane and over doubly periodic
cell domains.

\section{Non-Abelian vortex equations and existence theorems}
\setcounter{equation}{0}

We consider the non-Abelian Chern--Simons--Higgs vortex equations (\ref{mm1})--(\ref{mm2}).  With
\[
u =-\psi+\sum\limits_{i=1}^{M}\sum\limits_{s=1}^{n_i}\ln|z-z_{i, s}|^2,\quad  u_j=-\chi_j+\sum\limits_{s=1}^{n_j}\ln|z-z_{j, s}|^2,\quad j=1,\dots,M,
\]
we see that
 \[\re^u=\prod\limits_{i=1}^M\prod\limits_{s=1}^{n_i}|z-z_{i, s}|^2\re^{-\psi},\quad \re^{u_j}=\prod\limits_{s=1}^{n_j}|z-z_{j, s}|^2\re^{-\chi_j}, \quad j=1, \dots, M.\]
Then  the equations \eqref{mm1}--\eqref{mm2} become
  \ber
    \Delta u&=&\frac{\alpha^2}{M^2}\left(\sum\limits_{i=1}^M\left[\re^{u+u_i}+\re^{u-u_i}\right]-\xi\right)\left(\sum\limits_{j=1}^M\left[\re^{u+u_j}+\re^{u-u_j}\right]\right)\nm\\
    &&+\frac{\alpha\beta}{M}\sum\limits_{i=1}^M(\re^{u+u_i}-\re^{u-u_i})^2+4\pi\sum\limits_{i=1}^{M}\sum\limits_{s=1}^{n_i}\delta_{p_{i,s}}(x),\label{x3.1}\\
     \Delta u_j&=& \frac{\alpha\beta}{M}\left(\sum\limits_{i=1}^M\left[\re^{u+u_i}+\re^{u-u_i}\right]-\xi\right)(\re^{u+u_j}-\re^{u-u_j})\nm\\
     && +\beta^2(\re^{2u+2u_j}-\re^{2u-2u_j})+4\pi\sum\limits_{s=1}^{n_j}\delta_{p_{j,s}}(x),\quad
    j=1, \dots, M,\label{x3.2}
   \eer
 where  we set $\alpha=\frac\pi\kappa, \beta=\frac\pi\mu$, $p_{i, s}=z_{i, s}, i=1, \dots, M.$  When $M=1$ such equations were first obtained in
 \cite{gud1}.
 It will be convenient to take the rescaled parameters and translated variables
   \[\alpha\frac{\xi}{2M}\mapsto\alpha, \quad \beta\frac{\xi}{2M}\mapsto\beta,\quad u\mapsto u+\ln\frac{\xi}{2M}, \quad u_j\mapsto u_j, \,\,j=1, \dots, M. \]
Then the equations (\ref{x3.1}) and (\ref{x3.2}) are  renormalized into the form
   \ber
    \Delta u&=&\frac{\alpha^2}{M^2}\left(\sum\limits_{i=1}^M\left[\re^{u+u_i}+\re^{u-u_i}-2\right]\right)\left(\sum\limits_{j=1}^M\left[\re^{u+u_j}+\re^{u-u_j}\right]\right)\nm\\
    &&+\frac{\alpha\beta}{M}\sum\limits_{i=1}^M\left(\re^{u+u_i}-\re^{u-u_i}\right)^2
     +4\pi\sum\limits_{i=1}^{M}\sum\limits_{s=1}^{n_i}\delta_{p_{i,s}}(x),\label{e1}\\
     \Delta u_j&=& \frac{\alpha\beta}{M}\left(\sum\limits_{i=1}^M\left[\re^{u+u_i}+\re^{u-u_i}-2\right]\right)\left(\re^{u+u_j}-\re^{u-u_j}\right)\nm\\
     &&+\beta^2\left(\re^{2u+2u_j}-\re^{2u-2u_j}\right)
      +4\pi\sum\limits_{s=1}^{n_j}\delta_{p_{j,s}}(x),  \quad j=1, \dots, M.\label{e2}
   \eer

It will be interesting at this spot to compare the classical Abelian Chern--Simons--Higgs vortex equation (\ref{x7}) with the above system of non-Abelian vortex
equations, (\ref{e1}) and (\ref{e2}), in the Gudnason model \cite{gud1,gud2}, which will be our focus in the present work.

We will consider the equations \eqref{e1}--\eqref{e2} in two cases. In the first case we study the problem \eqref{e1}--\eqref{e2} over
the full plane $\mathbb{R}^2$  with  the topological boundary conditions
 \be
 u\to 0,\quad u_j\to 0   \quad \text{as }\quad |x|\to \infty, \quad j=1, \dots, M,\label{e3}
 \ee
realizing the asymptotic vacuum state with completely broken symmetry \cite{gud1,gud2}.

 We  have the following existence theorem.

\begin{theorem}\label{th1}
For any  sets of  points
\ber
Z_i=\{ p_{i,1}, \dots, p_{i,n_i}\}\subset \mathbb{R}^2, \quad i=1,\dots,M,
\eer
and the parameters   $\alpha, \beta>0,\, M\ge1$,  the system of
nonlinear elliptic equations \eqref{e1}--\eqref{e2} subject to the boundary condition \eqref{e3} admits a  solution over
$\mathbb{R}^2$ which possesses the quantized integrals
\ber
&&\frac{\alpha^2}{M^2}\int_{\mathbb{R}^2}\left(\sum\limits_{i=1}^M\left[\re^{u+u_i}+\re^{u-u_i}-2\right]\right)\left(\sum\limits_{j=1}^M\left[\re^{u+u_j}+\re^{u-u_j}\right]\right)\,{\rm\mbox{\rm d}} x\nm\\
    &&+\frac{\alpha\beta}{M}\sum\limits_{i=1}^M\int_{\mathbb{R}^2}\left(\re^{u+u_i}-\re^{u-u_i}\right)^2\,{\rm\mbox{\rm d}} x
     =-4\pi\sum\limits_{i=1}^{M}n_i,\label{q1}\\
&&\frac{\alpha\beta}{M}\int_{\mathbb{R}^2}\left(\sum\limits_{j=1}^M\left[\re^{u+u_j}+\re^{u-u_j}-2\right]\right)\left(\re^{u+u_i}-\re^{u-u_i}\right)\,{\rm\mbox{\rm d}} x\nm\\
     &&+\beta^2\int_{\mathbb{R}^2}\left(\re^{2u+2u_i}-\re^{2u-2u_i}\right)\,{\rm\mbox{\rm d}} x
      =-4\pi n_i,  \quad i=1, \dots, M.\label{q2}
\eer
Furthermore the boundary condition (\ref{e3}) is realized exponentially fast so that there hold the following asymptotic estimates near infinity:
\ber\label{asy}
u^2+\sum_{i=1}^M u_i^2 =\mbox{\rm O}(\re^{-m(1-\varepsilon)|x|}),\quad |\nabla u|^2+\sum_{i=1}^M |\nabla u_i|^2 =\mbox{\rm O}(\re^{-m(1-\varepsilon)|x|}),
\eer
where $m=2\sqrt{2}\min\{\alpha,\beta\}$ and $\varepsilon\in (0,1)$ is an arbitrarily small parameter.
\end{theorem}

In the second case we consider the equations \eqref{e1}--\eqref{e2} over a doubly periodic domain $\Omega$ with $M=1$. That is, in this
case we study a $2\times2$ version of  \eqref{e1}--\eqref{e2}.  For convenience we rewrite the system \eqref{e1}--\eqref{e2} with $M=1$ as follows
\ber
 \Delta U&=&\alpha^2\big(\re^{U+V}+\re^{U-V}\big)\left(\re^{U+V}+\re^{U-V}-2\right)
  +\alpha\beta\big(\re^{U+V}-\re^{U-V}\big)^2+4\pi\sum\limits_{j=1}^n\delta_{p_j}(x),\label{a1}\\
 \Delta V&=&\alpha\beta\big(\re^{U+V}-\re^{U-V}\big)\left(\re^{U+V}+\re^{U-V}-2\right)
  +\beta^2\big(\re^{2U+2V}-\re^{2U-2V}\big)+4\pi\sum\limits_{j=1}^n\delta_{p_j}(x).\label{a2}
 \eer

 We have the following existence results.

\begin{theorem}\label{thb1}
  Let $\Omega$ be a doubly periodic domain in $\mathbb{R}^2$ and $p_1, \dots,  p_n \in \Omega$ which need not to be distinct with repeated $p$'s counting for multiplicities.
Assume that  $\beta>\alpha>0$.
  \begin{enumerate}
    \item If
 the  equations \eqref{a1} and \eqref{a2} have a solution, then there holds
the condition
\be
8\pi n\leq \alpha\beta|\Omega|.
\ee
\item  Every solution $(U, V)$ of \eqref{a1} and \eqref{a2} satisfies
   \be
    \re^{U}< 1,\quad \re^{U+V}<1,\quad \re^{U-V}< 1.\label{b02}
    \ee
\item For any given  constant $\sigma>1$, assume
  \be
    \frac\beta\alpha<\sigma. \label{b03}
  \ee
   Then, there exist a positive constant $M_\sigma$ such that when $\alpha>M_\sigma$ the  equations
  \eqref{a1} and \eqref{a2} admit at least  two distinct   solutions  over $\Omega$, one of  which   satisfies the behavior
   \be
   \re^{U+V}\to 1, \quad\re^{U-V}\to 1, \quad as \quad \alpha \to +\infty \label{b04}
   \ee
 pointwise a.e. in $\Omega$. Furthermose, any solution $(U,V)$ of (\ref{a1}) and (\ref{a2}) possesses the quantized integrals
\be\label{xq1}
\alpha^2\int_{\Omega}\big(\re^{U+V}+\re^{U-V}\big)\left(\re^{U+V}+\re^{U-V}-2\right)\,\mbox{\rm d} x
  +\alpha\beta\int_{\Omega}\big(\re^{U+V}-\re^{U-V}\big)^2\,\mbox{\rm d} x=-4\pi n,
\ee
\be \label{xq2}
\alpha\beta\int_{\Omega}\big(\re^{U+V}-\re^{U-V}\big)\left(\re^{U+V}+\re^{U-V}-2\right)\,\mbox{\rm d} x
  +\beta^2\int_{\Omega}\big(\re^{2U+2V}-\re^{2U-2V}\big)\,\mbox{\rm d} x=-4\pi n.
\ee
\end{enumerate}
\end{theorem}

In Section 4, we establish Theorem \ref{th1} using a direct minimization method which extends the techniques in \cite{jata,tau1,Ycs,Ybook}.
In Section 5, we prove Theorem \ref{thb1} by utilizing and extending an inequality-constrained variational method originally developed in \cite{caya1} and further developed in  \cite{taran96, nota1,nota2},  and subsequently in \cite{nota} for a context that relates more to the situation considered here.   In this respect, see also \cite{hata}.

\section{Proof of Theorem 3.1}
\setcounter{equation}{0}

In this section we prove the existence of topological solutions for the equations \eqref{e1}--\eqref{e2}.

Choosing  the background functions
  \ber \label{x4.1}
u_i^0(x)=-\sum\limits_{s=1}^{n_i}\ln\big(1+\lambda|x-p_{i, s}|^{-2}\big), \quad \lambda>0,\quad i=1, \dots, M,\eer
which satisfy
  \be\label{x4.2}
 \Delta u_i^0=-h_i+4\pi\sum\limits_{s=1}^{n_i}\delta_{p_{i,s}}, \quad h_i(x)=4\lambda\sum\limits_{s=1}^{n_i}\frac{1}{(\lambda+|x-p_{i,s}|^2)^2},\ee
we see that   the new variables
 \be
  u=\sum\limits_{i=1}^M u_i^0+f, \quad u_j=u_j^0+f_j, \quad j=1, \dots, M, \label{e4}
 \ee
 allow us to recast the equations \eqref{e1}--\eqref{e2} into
    \ber
    \Delta f&=&\frac{\alpha^2}{M^2}\left(\sum\limits_{i=1}^M\left[\re^{\sum\limits_{k=1}^M u_k^0+u_i^0+f+f_i}+\re^{\sum\limits_{k=1}^M u_k^0-u_i^0+f-f_i}-2\right]\right)\left(\sum\limits_{j=1}^M\left[\re^{\sum\limits_{k=1}^M u_k^0+u_j^0+f+f_j}+\re^{\sum\limits_{k=1}^M u_k^0-u_j^0+f-f_j}\right]\right)
    \nm\\  &&+\frac{\alpha\beta}{M}\sum\limits_{i=1}^M\left(\re^{\sum\limits_{k=1}^M u_k^0+u_i^0+f+f_i}-\re^{\sum\limits_{k=1}^M u_k^0-u_i^0+f-f_i}\right)^2+\sum\limits_{i=1}^Mh_i,\label{e5}\\
     \Delta f_i&=& \frac{\alpha\beta}{M}\left(\sum\limits_{j=1}^M\left[\re^{\sum\limits_{k=1}^M u_k^0+u_j^0+f+f_j}+\re^{\sum\limits_{k=1}^M u_k^0-u_j^0+f-f_j}-2\right]\right)\left(\re^{\sum\limits_{k=1}^M u_k^0+u_i^0+f+f_i}-\re^{\sum\limits_{k=1}^M u_k^0-u_i^0+f-f_i}\right)
      \nm\\  &&+\beta^2\left(\re^{2\sum\limits_{k=1}^M u_k^0+2u_i^0+2f+2f_i}-\re^{2\sum\limits_{k=1}^M u_k^0-2u_i^0+2f-2f_i}\right)+h_i, \quad i=1, \dots, M.\label{e6}
   \eer

The topological boundary condition \eqref{e3}  becomes
 \be
 f\to 0,\quad f_i\to 0   \quad \text{as }\quad |x|\to \infty, \quad i=1,\dots, M.\label{e7}
 \ee

It can be checked that the  equations (\ref{e5}) and (\ref{e6}) are the   Euler--Lagrange equations of the   functional
   \ber
    &&I(f, f_1, \dots, f_M)\nm\\
    &&=\itr \ud x\left\{\frac{M}{\alpha}|\nabla f|^2+\frac1\beta\sum\limits_{i=1}^M|\nabla f_i|^2
    +\sum\limits_{i=1}^M\left(\frac\alpha M\left[\re^{\sum\limits_{k=1}^M u_k^0+u_i^0+f+f_i}+\re^{\sum\limits_{k=1}^M u_k^0-u_i^0+f-f_i}-2\right]^2\right.\right.\nm\\
    &&\left.\left.+\beta\left[\re^{\sum\limits_{k=1}^M u_k^0+u_i^0+f+f_i}-\re^{\sum\limits_{k=1}^M u_k^0-u_i^0+f-f_i}\right]^2\right)+\frac{2M}{\alpha}\sum\limits_{i=1}^Mfh_i
    +\frac{2}{\beta}\sum\limits_{i=1}^Mf_ih_i\right\}.\label{e8}
    \eer

We consider the functional $I$ over $\wot$.  Here and in what follows we use  $\wot$ to denote the usual Sobolev space
of scalar-valued or vector-valued functions.  It is not difficult to see that the functional $I$ is continuous, differentiable and lower
semi-continuous on $\wot$. The important thing is that we can show  that the functional $I$ is coercive and bounded from below  over $\wot$, which will be carried out later. Then we can conclude  that the functional $I$ admits a critical point   $(f,f_1, \dots, f_M)\in\wot$, which is  a weak solution to the equations
 \eqref{e5}--\eqref{e6}.  By the following inequality
  \[\|\re^w-1\|^2_2\le C\exp\left(C\|w\|_{\wot}^2\right),\quad \forall\, w\in \wot,\]
we see that the right hand side of the equations \eqref{e5}--\eqref{e6} belongs to $L^2(\mathbb{R}^2)$. Then  using  elliptic $L^2$-estimates and a bootstrap argument, we find that the
solution $(f,f_1, \dots, f_M)$ is smooth. In particular, $(f,f_1, \dots, f_M)$ lies in $W^{2,2}(\mathbb{R}^2)$ which ensures that $(f,f_1, \dots, f_M)$ satisfies the
boundary condition \eqref{e7}. Then, in view of \eqref{e4}, we see that the problem consisting of the equations
\eqref{e1}  and \eqref{e2} admits a solution $(u, u_1,\dots, u_M)$ satisfying the boundary condition \eqref{e3}.

We now establish the exponential decay estimates for the solution.

We first note that
\ber
\re^{u+u_i}+\re^{u-u_i}-2&=&(\re^{u+u_i}-1)+(\re^{u-u_i}-1)=\re^{\xi_i'}(u+u_i)+\re^{\xi_i''}(u-u_i)\nm\\
&=&(\re^{\xi_i'}+\re^{\xi_i'})u +(\re^{\xi_i'}-\re^{\xi_i''})u_i,\label{xx1}\\
\re^{u+u_i}-\re^{u-u_i}&=&2\re^{\xi_i}u_i,\label{xx2}
\eer
where $\xi_i'$, $\xi_i''$, and $\xi_i$ are between $0$ and $u+u_i$, $0$ and $u-u_i$, and $u+u_i$ and $u-u_i$, respectively, $i=1,\dots,M$.

By virtue of (\ref{xx1}) and (\ref{xx2}),  we see that the equations  \eqref{e1}--\eqref{e2} may be
rewritten as
  \ber
\Delta u&=&\frac{\alpha^2}{M^2}\left(\sum_{j=1}^M\left[\re^{u+u_j}+\re^{u-u_j}\right]\right)\left(\sum_{i=1}^M\left[(\re^{\xi_i'}+\re^{\xi_i'})u +(\re^{\xi_i'}-\re^{\xi_i''})u_i\right]\right)\nm\\
&&+\frac{4\alpha\beta}M\sum_{i=1}^M\re^{2\xi_i}u_i^2,\label{xe1}\\
\Delta u_i&=&\frac{2\alpha\beta}M\re^{\xi_i} \left(\sum_{j=1}^M\left[(\re^{\xi_j'}+\re^{\xi_j'})u +(\re^{\xi_j'}-\re^{\xi_j''})u_j\right]\right)u_i\nm\\
&&+2\beta^2\left(\re^{u+u_i}+\re^{u-u_i}\right)\re^{\xi_i} u_i,\quad i=1,\dots,M,\label{xe2}
\eer
where $|x|>R$ and $R>0$ is taken to be sufficiently large so that $R>|p_{i,s}|$ for $s=1,\dots,n_i$ and $i=1,\dots,M$.

To proceed further, we set
\ber \label{xU}
U=u^2+\sum_{i=1}^M u_i^2.
\eer
Then we can compute for $|x|>R$ the result
\ber\label{xe3}
\Delta U&\geq& 2u\Delta u +2\sum_{i=1}^M u_i\Delta u_i\nm\\
&\geq&\frac{2\alpha^2}{M^2}\left(\sum_{j=1}^M\left[\re^{u+u_j}+\re^{u-u_j}\right]\right)\left(\sum_{i=1}^M\left[\re^{\xi_i'}+\re^{\xi_i'}\right]\right)u^2 +4\beta^2
\sum_{i=1}^M \left(\re^{u+u_i}+\re^{u-u_i}\right)\re^{\xi_i} u_i^2\nm\\
&&-\frac{2\alpha^2}{M^2}\left(\sum_{j=1}^M\left[\re^{u+u_j}+\re^{u-u_j}\right]\right)\left(\sum_{i=1}^M |\re^{\xi_i'}-\re^{\xi_i''}||u_i|\right)|u|\nm\\
&& -\frac{8\alpha\beta}M\sum_{i=1}^M\re^{2\xi_i}u_i^2|u|-\frac{4\alpha\beta}M\sum_{i=1}^n\re^{\xi_i} \left(\sum_{j=1}^M\left[(\re^{\xi_j'}+\re^{\xi_j'})|u| +|\re^{\xi_j'}-\re^{\xi_j''}||u_j|\right]\right)u_i^2.
\eer
Applying (\ref{e3}) and the Schwartz inequality in (\ref{xe3}),  we see that for any arbitrarily small $\varepsilon\in(0,1)$ there is $R_\varepsilon>R$ such that $U$ satisfies the
elliptic inequality
\ber \label{xe4}
\Delta U\geq 8(\min\{\alpha,\beta\})^2\left(1-\frac\varepsilon2\right) U,\quad |x|\geq R_\varepsilon.
\eer
Applying a comparison function argument to (\ref{xe4}) and using the boundary property $U=0$ at infinity, we can find a sufficient large constant $C(\varepsilon)>0$ such that
\ber\label{xe5}
U(x)\leq C(\varepsilon)\re^{-2\sqrt{2}\min\{\alpha,\beta\}(1-\varepsilon)|x|},\quad |x|\geq R_\varepsilon.
\eer

Next, using $\partial$ to denote one of the two partial derivatives, $\partial_1$ and $\partial_2$, we obtain from (\ref{e1}) and (\ref{e2}) the results
\ber
\Delta(\partial u)&=&\frac{\alpha^2}{M^2}\left(\sum_{i=1}^M\left[\re^{u+u_i}+\re^{u-u_i}\right]\right)^2\pa u
+\frac{\alpha^2}{M^2}\left(\sum_{i=1}^M\left[\re^{u+u_i}+\re^{u-u_i}\right]\right)\left(\sum_{i=1}^M\left[\re^{u+u_i}-\re^{u-u_i}\right]\right)\pa u_i\nm\\
&&+\frac{\alpha^2}{M^2}\left(\sum_{i=1}^M\left[\re^{u+u_i}+\re^{u-u_i}-2\right]\right)\left(\sum_{j=1}^M\left[\re^{u+u_j}+\re^{u-u_j}\right]\pa u+\sum_{j=1}^M\left[\re^{u+u_j}-\re^{u-u_j}\right]
\pa u_j\right)\nm\\
&&+\frac{2\alpha\beta}M\sum_{i=1}^M\left(\re^{u+u_i}-\re^{u-u_i}\right)\left(\left[\re^{u+u_i}-\re^{u-u_i}\right]\pa u+\left[\re^{u+u_i}+\re^{u-u_i}\right]\pa u_i\right),\label{ep1}\\
\Delta(\pa u_i)&=&2\beta^2\left(\re^{2u+2u_i}+\re^{2u-2u_i}\right)(\pa u_i)+2\beta^2\left(\re^{2u+2u_i}-\re^{2u-2u_i}\right)(\pa u)\nm\\
&&+\frac{\alpha\beta}M\left(\sum_{j=1}^M\left[\re^{u+u_j}+\re^{u-u_j}-2\right]\right)\left(\left[\re^{u+u_i}-\re^{u-u_i}\right]\pa u+\left[\re^{u+u_i}+\re^{u-u_i}\right]\pa u_i\right)\nm\\
&&+\frac{\alpha\beta}M\left(\re^{u+u_i}-\re^{u-u_i}\right)\left(\sum_{j=1}^M\left[\re^{u+u_j}+\re^{u-u_j}\right]\pa u+\sum_{j=1}^M\left[\re^{u+u_j}-\re^{u-u_j}\right]
\pa u_j\right).\label{ep2}
\eer
Applying the $L^2$-estimate and the fact that $u,u_1,\dots,u_M\in W^{2,2}$ outside $B_R=\{x\in\mathbb{R}^2\,|\, |x|>R\}$, we see in view of the above equations that $\pa u,\pa u_1,\dots,\pa u_M\in W^{2,2}$ outside
$B_R$ as well. Therefore
\be
\pa u,\pa u_1,\dots,\pa u_M \to0\quad \mbox{as}\quad |x|\to\infty.
\ee

As before, we set
\be
V=(\pa u)^2+\sum_{i=1}^M(\pa u_i)^2.
\ee
Similar to the case with the function $U$ defined in (\ref{xU}), we may apply the Schwartz inequality and use the equations (\ref{ep1}) and (\ref{ep2}) to obtain the elliptic inequality
\be
\Delta V\geq8(\min\{\alpha,\beta\})^2\left(1-\frac\varepsilon2\right) V,\quad |x|\geq R_\varepsilon.
\ee
Thus $V$ enjoys the same exponential decay estimate as $U$ as stated in (\ref{xe5}).

Now we only need to prove the coerciveness and a bound from below  for  the functional $I(f,  f_1, \dots, f_M)$ over $\wot$.

   Using the elementary inequality
   \[ \tilde{\alpha}(a+b)^2+\tilde{\beta}(a-b)^2\ge 2\min\{\tilde{\alpha}, \tilde{\beta}\}(a^2+b^2),
   \quad \forall\,\tilde{\alpha},\tilde{\beta}>0,\, \forall \,a, b\in \mathbb{R},\]
 with
  \[a=\re^{\sum\limits_{k=1}^M u_k^0+u_i^0+f+f_i}-1, \quad b=\re^{\sum\limits_{k=1}^M u_k^0-u_i^0+f-f_i}-1, \]
   we see that the functional $I(f, f_1, \dots, f_M)$ over $\wot$ satisfies
  \ber
    &&I(f, f_1, \dots, f_M)\nm\\
    &&\ge\itr \ud x\left\{\frac{M}{\alpha}|\nabla f|^2+\frac1\beta\sum\limits_{i=1}^M|\nabla f_i|^2
    +2\min\left\{\frac\alpha M, \beta\right\}\sum\limits_{i=1}^M\left(\re^{\sum\limits_{k=1}^M u_k^0+u_i^0+f+f_i}-1\right)^2\right.\nm\\
    &&\left.\quad+2\min\left\{\frac\alpha M, \beta\right\}\sum\limits_{i=1}^M\left(\re^{\sum\limits_{k=1}^M u_k^0-u_i^0+f-f_i}-1\right)^2+\frac{2M}{\alpha}\sum\limits_{i=1}^Mfh_i
    +\frac{2}{\beta}\sum\limits_{i=1}^Mf_ih_i\right\}.\label{e9}
    \eer

To  proceed we need the  following lemma in \cite{wangr}.
\begin{lemma}
 The function $h_i$ belongs to $ L^2(\mathbb{R}^2)$  with
  \be
   \|h_i\|_2\le \frac{C}{\sqrt{\lambda}},\label{e10}
  \ee
for some positive constant $C$ independent of $\lambda$ and $\re^{u_i^0}-1\in L^p(\mathbb{R}^2)$ for any  $p\ge2,  i=1, \dots, M$.
\end{lemma}
Here and in what  follows we use $C$ to denote a positive constant which may take different values at different places.

Using the H\"{o}lder inequality  and \eqref{e10}, we have
  \ber
  \itr \left(\frac{2M}{\alpha}\sum\limits_{i=1}^Mfh_i+\frac{2}{\beta}\sum\limits_{i=1}^Mf_ih_i\right)\ud x&\ge&
  -\frac{2M}{\alpha}\sum\limits_{i=1}^M\|f\|_2\|h_i\|_2-\frac{2}{\beta}\sum\limits_{i=1}^M\|f_i\|_2\|h_i\|_2\nm\\
  &\ge&-\frac{C}{\sqrt{\lambda}}\left(\|f\|_2+\sum\limits_{i=1}^M\|f_i\|_2\right).\label{e11}
  \eer

   In what follows  we need to control the $L^2$-norm of $f$ and $f_j\,(j=1, \dots, M)$ by the  positive terms in  \eqref{e9}.
 Now we deal with the third term on the right-hand side of \eqref{e9}.  Since
  \[\re^{u_i^0}-1\in L^2(\mathbb{R}^2),\quad i=1, \dots, M, \]
 it is easy to check that
  \be
  \re^{\sum\limits_{k=1}^M u_k^0+u_i^0}-1\in L^2(\mathbb{R}^2), \quad \re^{\sum\limits_{k=1}^M u_k^0-u_i^0}-1\in L^2(\mathbb{R}^2), \quad i=1, \dots, M.\label{e15a}
  \ee
  To proceed further, we need   the inequality
\be
 |\re^t-1|\ge \frac{|t|}{1+|t|}, \quad \forall\, t\in \mathbb{R}, \label{e12}
 \ee
 which follows directly from the elementary inequalities
$\re^t-1\ge t, \forall\, t\ge0$,
 and
 $1-\re^{-t}\ge \frac{t}{1+t}, \forall\,t\ge0.$

 With \eqref{e15a} and \eqref{e12}, we have
 \ber
 \itr\left(\re^{\sum\limits_{k=1}^M u_k^0+u_i^0+f+f_i}-1\right)^2&=&\itr\left(\re^{\sum\limits_{k=1}^M u_k^0+u_i^0}\left[\re^{f+f_i}-1\right]+\re^{\sum\limits_{k=1}^M u_k^0+u_i^0}-1\right)^2\ud x\nm\\
 &\ge& \frac12\itr \re^{2\sum\limits_{k=1}^M u_k^0+2u_i^0}\left(\re^{f+f_i}-1\right)^2-\itr\left(\re^{\sum\limits_{k=1}^M u_k^0+u_i^0}-1\right)^2\ud x\nm\\
  &\ge& \frac12\itr \re^{2\sum\limits_{k=1}^M u_k^0+2u_i^0} \frac{|f+f_i|^2}{(1+|f+f_i|)^2}\ud x-C, \quad i=1, \dots, M. \label{e17}
 \eer

By the definition of $u_i^0$, we see that $ \re^{2\sum\limits_{k=1}^M u_k^0+2u_i^0}$ satisfies    $0\le \re^{2\sum\limits_{k=1}^M u_k^0+2u_i^0}<1$,  vanishes  at
the vortex point $p_{i,s},  (s=1, \dots, n_i,\, i=1, \dots, M)$, and approaches  $1$ at infinity.  As  in \cite{jata,tau1}, we decompose
 $\mathbb{R}^2$ as follows,
  \be
   \mathbb{R}^2=\Omega_1^i\cup\Omega_2^i, \quad i=1, \dots, M,\label{e18}
   \ee
where \[ \quad \Omega_1^i=\left\{x\in\mathbb{R}^2\bigg| \, \re^{2\sum\limits_{k=1}^M u_k^0+2u_i^0}\le \frac12 \right\}, \quad \Omega_2^i=\left\{x\in\mathbb{R}^2\bigg| \, \re^{2\sum\limits_{k=1}^M u_k^0+2u_i^0}\ge \frac12 \right\},\quad i=1, \dots, M. \]

 To deal with the right-hand side of \eqref{e17}, we need the inverse H\"{o}lder inequality (cf. \cite{wangr}):
\begin {lemma}
 For any measurable functions $g_1, g_2$ on the domain $\Omega$,  there holds the inequality
 \be
   \ito|g_1g_2|\ud x \ge \left(\ito |g_1|^q\ud x\right)^{\frac1q}\left(\ito |g_2|^{q'}\ud x\right)^{\frac{1}{q'}},\label{e19}
  \ee
 where $q,q'\in \mathbb{R}$, $0<q<1$, $q'<0$ and $\frac1q+\frac{1}{q'}=1.$
\end {lemma}

 On $\Omega_1^i$, we have  $0\le\re^{2\sum\limits_{k=1}^M u_k^0+2u_i^0}\le\frac12$   and $\re^{2\sum\limits_{k=1}^M u_k^0+2u_i^0}$ tends to $0$ at most  at order $4\sum\limits_{k=1}^Mn_k+4n_i$  near the
 vortex points. Then, by  taking  $q'_i$ satisfying
\[
-\frac{1}{2\sum\limits_{k=1}^Mn_k+2n_i}<q'_i<0,
\]
 we see that the integrals
  \[\int_{\Omega_1^i}\re^{2q'_i\sum\limits_{k=1}^M u_k^0+2q'_iu_i^0}\ud x, \quad i=1, \dots, M\] exist.

 Using the inverse H\"{o}lder inequality \eqref{e19}, we can get
  \ber
  \int_{\Omega_1^i} \re^{2\sum\limits_{k=1}^M u_k^0+2u_i^0} \frac{|f+f_i|^2}{(1+|f+f_i|)^2}\ud x
  &\ge&\left(\int_{\Omega_1^i}\frac{|f+f_i|^{2q_i}}{(1+|f+f_i|)^{2q_i}}\ud x\right)^{\frac{1}{q_i}}
   \left(\int_{\Omega_1^i} \re^{2q'_i\sum\limits_{k=1}^M u_k^0+2q'_iu_i^0}\ud x\right)^{\frac{1}{q'_i}}\nm\\
   &\ge&C \left(\int_{\Omega_1^i}\frac{|f+f_i|^{2q_i}}{(1+|f+f_i|)^{2q_i}}\ud x\right)^{\frac{1}{q_i}},\quad i=1, \dots, M,\label{e20}
  \eer
where
\[0<q_i<\frac{1}{2\sum\limits_{k=1}^Mn_k+2n_i+1},\quad i=1, \dots, M.\]
Noting
\[0\le\frac{|f+f_i|}{1+|f+f_i|}<1, \quad i=1, \dots, M\]
and applying the Young inequality, we obtain
 \ber
  \left(\int_{\Omega_1^i}\frac{|f+f_i|^{2q_i}}{(1+|f+f_i|)^{2q_i}}\ud x\right)^{\frac{1}{q_i}}
  &\ge&\left(\int_{\Omega_1^i}\frac{|f+f_i|^2}{(1+|f+f_i|)^2}\ud x\right)^{\frac{1}{q_i}} \nm\\
  &\ge& C\int_{\Omega_1^i}\frac{|f+f_i|^2}{(1+|f+f_i|)^2}\ud x-C, \quad i=1, \dots, M.\label{e21}
 \eer
 Combining \eqref{e20} and \eqref{e21}, we have
 \be
  \int_{\Omega_1^i} \re^{2\sum\limits_{k=1}^M u_k^0+2u_i^0} \frac{|f+f_i|^2}{(1+|f+f_i|)^2}\ud x\ge C\int_{\Omega_1}\frac{|f+f_i|^2}{(1+|f+f_i|)^2}\ud x-C,\quad i=1, \dots, M. \label{e22}
 \ee

 On the other hand, over $\Omega_2^i$, it is easy to get
  \be
  \int_{\Omega_2^i} \re^{2\sum\limits_{k=1}^M u_k^0+2u_i^0} \frac{|f+f_i|^2}{(1+|f+f_i|)^2}\ud x\ge \frac12\int_{\Omega_2^i}\frac{|f+f_i|^2}{(1+|f+f_i|)^2}\ud x, \quad i=1, \dots, M.\label{e23}
  \ee
 Hence, from \eqref{e17}, \eqref{e22}, and \eqref{e23}, we infer that
  \be
  \itr\left(\re^{\sum\limits_{k=1}^M u_k^0+u_i^0+f+f_i}-1\right)^2\ge C\itr\frac{|f+f_i|^2}{(1+|f+f_i|)^2}\ud x-C, \quad i=1, \dots, M.\label{e24}
  \ee

 Now repeating  the procedure in getting \eqref{e24},   we have
  \be
  \itr\left(\re^{\sum\limits_{k=1}^M u_k^0-u_i^0+f-f_i}-1\right)^2\ge C\itr\frac{|f-f_i|^2}{(1+|f-f_i|)^2}\ud x-C, \quad i=1, \dots, M.\label{e24a}
  \ee

 To proceed further, we invoke the following standard interpolation inequality over $\mathbb{R}^2$:
 \be
 \itr w^4\ud x\le 2\itr w^2\ud x\itr |\nabla w|^2\ud x, \quad \forall\, w\in \wot. \label{e14}
 \ee
 Using \eqref{e14}, we obtain
  \ber
   &&\left(\itr|f+f_i|^2\ud x\right)^2\nm\\
   &&=\left(\itr\frac{|f+f_i|}{1+|f+f_i|}[1+|f+f_i|]|f+f_i|\ud  x\right)^2\nm\\
    &&\le\itr\frac{|f+f_i|^2}{(1+|f+f_i|)^2}\ud x\itr\big(|f+f_i|+|f+f_i|^2\big)^2\ud x\nm\\
    &&\le4\itr\frac{|f+f_i|^2}{(1+|f+f_i|)^2}\ud x\itr|f+f_i|^2\ud x\left(\itr|\nabla(f+f_i)|^2+1\right)\nm\\
    &&\le\frac12 \left(\itr|f+f_i|^2\ud x\right)^2+C\left(\left[\itr\frac{|f+f_i|^2}{[1+|f+f_i|]^2}\ud x\right]^4+\left[\itr|\nabla[f+f_i]|^2\ud x\right]^4+1\right), \label{e15}
  \eer
  which  implies
   \ber
   \|f+f_i\|_2&\le& C\left(\itr \frac{|f+f_i|^2}{[1+|f+f_i|]^2}\ud x+\itr|\nabla[f+f_i]|^2\ud x+1\right)\nm\\
   &\le& C\left(\itr \frac{|f+f_i|^2}{[1+|f+f_i|]^2}\ud x+\itr\left[|\nabla f|^2+|\nabla f_i|^2\right]\ud x+1\right), \quad i=1, \dots, M. \label{e16}
   \eer

 Similarly, we have
  \be
   \|f-f_i\|_2\le C\left(\itr \frac{|f-f_i|^2}{(1+|f-f_i|)^2}\ud x+\itr\left[|\nabla f|^2+|\nabla f_i|^2\right]\ud x+1\right), \quad i=1, \dots, M.   \label{e25}
  \ee

 Then, in view of \eqref{e16}, \eqref{e25},  and  the following  simple  inequality
   \[\|f\|_2+\|f_i\|_2\le 2(\|f+f_i\|_2+\|f-f_i\|_2), \quad i=1, \dots, M, \]
    we see that
 \ber
  &&\|f\|_2+\sum\limits_{i=1}^M\|f_i\|_2\nm\\
  &&\le  C\left(\itr|\nabla f|^2\ud x+\sum\limits_{i=1}^M\itr\left[|\nabla f_i|^2+\frac{|f+f_i|^2}{(1+|f+f_i|)^2}+ \frac{|f-f_i|^2}{(1+|f-f_i|)^2}\right]\ud x+1\right).  \label{e26}
 \eer

From \eqref{e9}, \eqref{e11}, \eqref{e24} and \eqref{e24a}, we conclude that
 \ber
 &&I(f, f_1, \dots, f_M)\nm\\
 &&\ge C\left( \itr|\nabla f|^2\ud x+\sum\limits_{i=1}^M\itr\left[|\nabla f_i|^2+ \frac{|f+f_i|^2}{(1+|f+f_i|)^2}+\frac{|f-f_i|^2}{(1+|f-f_i|)^2}\right]\ud x\right)\nm\\
 &&\quad-\frac{C}{\sqrt{\lambda}}\left(\|f\|_2+\sum\limits_{i=1}^M\|f_i\|_2\right)-C.\label{e27}
 \eer

At this point,  combining \eqref{e26} and \eqref{e27} and taking $\lambda$ sufficiently large, we can get

\ber
   &&I(f, f_1,\dots, f_M)\nm\\
   &&\ge C\left( \itr|\nabla f|^2\ud x+\sum\limits_{i=1}^M\itr\left[|\nabla f_i|^2+ \frac{|f+f_i|^2}{(1+|f+f_i|)^2}+\frac{|f-f_i|^2}{(1+|f-f_i|)^2}\right]\ud x\right)-C.\label{e28}
\eer
 Then applying \eqref{e26} in the right hand side of \eqref{e28},
 we have
 \be
  I(f, f_1, \dots, f_M)\ge C\left(\|f\|_{\wot}+\sum\limits_{i=1}^M\|f_i\|_{\wot}\right)-C,
  \ee
which says that the functional $I(f, f_1,\dots, f_M)$ is coercive and bounded from below  over $\wot$. Therefore the existence of a critical point as a global minimizer of $I$ in $\wot$ follows immediately.

In order to establish the results regarding the quantized integrals (\ref{q1}) and (\ref{q2}), we note that the background functions $u_i^0$ ($i=1,\dots,M$) defined in
(\ref{x4.1}) obey the decay estimates
\be
|\nabla u_i^0(x)|=\mbox{O}(|x|^{-3})\quad\mbox{as }|x|\to\infty,\quad i=1,\dots,M.
\ee
On the other hand, since the solution $(u,u_1,\dots,u_M)$ of (\ref{e1})--(\ref{e2}) obtained decays at infinity according to (\ref{asy}), we see that $(f,f_1,\dots,f_M)$
set forth in (\ref{e4}) satisfies
\be \label{asy2}
|\nabla f(x)|+\sum_{i=1}^M |\nabla f_i(x)|=\mbox{O}(|x|^{-3})\quad\mbox{as }|x|\to\infty.
\ee
Using (\ref{asy2}) and the divergence theorem, we arrive at
\be \label{4.45}
\int_{\mathbb{R}^2}\Delta f\,\mbox{d} x=0,\quad\int_{\mathbb{R}^2}\Delta f_i\,\mbox{d} x=0,\quad i=1,\dots,M.
\ee
Moreover, integrating directly, we have
\be \label{4.46}
\int_{\mathbb{R}^2}h_i\,\mbox{d}x=4\pi n_i,\quad i=1,\dots,M.
\ee
Finally, integrating (\ref{e5}) and (\ref{e6}) over $\mathbb{R}^2$ and applying (\ref{4.45}) and (\ref{4.46}), we obtain the quantized integrals (\ref{q1}) and (\ref{q2}) stated
in the theorem.

The proof of Theorem \ref{th1} is now complete.

\section{Proof of Theorem 3.2}
\setcounter{equation}{0}\setcounter{lemma}{0}

In this section we establish the existence of solutions to \eqref{a1}--\eqref{a2} over a doubly periodic domain.  We will
make a variational formulation of the problem. Then we can carry out a constrained  minimization procedure to find the critical
points for the associated functional.  The key step is to find some inequality-type constraints, from which we can define a suitable
admissible set.   This procedure was  initiated in \cite{caya1} and refined  in \cite{nota1,nota2,nota} and \cite{hata}.

We first give {\em a priori}  estimates of the solutions  to  \eqref{a1}--\eqref{a2}.
 \begin{proposition} \label{pr1}
  Let $(U, V)$ be a solution of \eqref{a1}--\eqref{a2}. Then $U<0,   U+V<0, U-V<0$ throughout $\Omega$.
  \end{proposition}

{\bf Proof.} \quad Let $(U, V)$ be a solution of \eqref{a1} and \eqref{a2}.   Introduce a transformation $f=U+V, g=U-V$. From \eqref{a1} and \eqref{a2}, we conclude that $f$ and $g$
satisfy the  equations
  \ber
 \Delta f&=&(\alpha+\beta)^2\re^f\big(\re^f-1\big)+(\alpha-\beta)^2\re^g\big(\re^f-1\big)\nm\\
&&-(\beta^2-\alpha^2)\big(\re^f+\re^g\big)\big(\re^g-1\big)
  +8\pi\sum\limits_{j=1}^n\delta_{p_j},\label{b1'}\\
   \Delta g&=&(\alpha+\beta)^2\re^g\big(\re^g-1\big)+(\alpha-\beta)^2\re^f\big(\re^g-1\big)-(\beta^2-\alpha^2)\big(\re^f+\re^g\big)\big(\re^f-1\big).\label{b2'}
  \eer
  We first show that $U<0$ in $\Omega$. From \eqref{a1}, we see that
   \berr
    \Delta U&\ge& \alpha^2\re^U\big(\re^V+\re^{-V}-2\re^{-U}\big)\re^U\big(\re^V+\re^{-V}\big)+4\pi\sum\limits_{j=1}^n\delta_{p_j}\\
     &\ge&2\alpha^2\re^U\big(\re^V+\re^{-V}\big)\big(\re^U-1\big)+4\pi\sum\limits_{j=1}^n\delta_{p_j}.
   \eerr
Then, by maximum principle, we have $U<0$  throughout $\Omega$.

To prove $U+V<0,$  we argue by contradiction.  Assume that there exists a point
 $\tilde{x}\in\Omega$ such that
   \[f(\tilde{x})=\max\limits_{x\in\Omega}f(x)\ge0.\]
   From the equation \eqref{b1'}, we  have  $g(\tilde{x})\ge0$.
   Then  we obtain
   $U(\tilde{x})=\frac12(f(\tilde{x})+g(\tilde{x}))\ge0$, which
   contradicts the conlusion $U<0$ in $\Omega$. Therefore,  we have $f<0$
   in $\Omega$.

   Similarly, if there is a point $\tilde{x}\in \Omega$ such that
  \[g(\tilde{x})=\max\limits_{x\in\Omega}g(x)\ge0,\]
then by the equation \eqref{b2'}, we see that $f(\tilde{x})\ge0$, which
again leads to a contradiction. Hence the conclusion follows.

By Proposition \ref{pr1}, the second part of Theorem \ref{thb1}
follows.

 Let $u_0$ be the  unique solution of the following problem (see
\cite{aubi})
  \berr
 \Delta u_0=-\frac{8\pi n}{|\Omega|}+8\pi\sum\limits_{j=1}^n\delta_{p_j} \quad \text{on } \Omega;\quad
 \ito u_0\ud x=0.
  \eerr

For convenience,  we introduce the following new variables:
  \be U=\frac{u_0}{2}+\frac{u+v}{2},\quad   V=\frac{u_0}{2}+\frac{u-v}{2}, \label{b3}
  \ee
which reduce the equations   \eqref{a1}--\eqref{a2}   into  the
  form
  \ber
   \Delta\frac{u+v}{2}&=&\alpha^2\big(\re^{u_0+u}+\re^{v}\big)\left(\re^{u_0+u}+\re^{v}-2\right)
  +\alpha\beta\big(\re^{u_0+u}-\re^{v}\big)^2+\frac{4\pi n}{|\Omega|},\label{b4}\\
 \Delta\frac{u-v}{2} &=&\alpha\beta\big(\re^{u_0+u}-\re^{v}\big)\left(\re^{u_0+u}+\re^{v}-2\right)
  +\beta^2\big(\re^{2u_0+2u}-\re^{2v}\big)+\frac{4\pi n}{|\Omega|}.\label{b5}
  \eer

  To make a variational reformulation of the problem, we rewrite \eqref{b4} and \eqref{b5}  equivalently as
  \ber
   \frac12\left(\frac1\alpha+\frac1\beta\right)\Delta u+\frac12\left(\frac1\alpha-\frac1\beta\right)\Delta v&=&2\re^{u_0+u}\left([\alpha+\beta]\left[\re^{u_0+u}-1\right]+[\alpha-\beta]\left[\re^v-1\right]\right)
   \nm\\&&+\frac{4\pi n}{|\Omega|}\left(\frac1\alpha+\frac1\beta\right),\label{b4a}\\
   \frac12\left(\frac1\alpha-\frac1\beta\right)\Delta u+\frac12\left(\frac1\alpha+\frac1\beta\right)\Delta v&=&2\re^v\left([\alpha-\beta]\left[\re^{u_0+u}-1\right]+[\alpha+\beta]\left[\re^v-1\right]\right)
   \nm\\&&+\frac{4\pi n}{|\Omega|}\left(\frac1\alpha-\frac1\beta\right).\label{b5a}
  \eer

  Therefore, in the sequel we only need to solve
  \eqref{b4a} and \eqref{b5a}.

We will work on the space $W^{1, 2}(\Omega)\times W^{1, 2}(\Omega)$, where $W^{1, 2}(\Omega)$ denotes the set of $\Omega$-periodic $L^2$-
functions whose  derivatives are also in $L^2(\Omega)$. We denote the usual  norm on $W^{1, 2}(\Omega)$ by $\|\cdot\|$ as given by
 \[
\|w\|^2=\|w\|_2^2+\|\nabla w\|_2^2=\ito w^2\ud x+\ito |\nabla w|^2\ud x.
\]

 It is easy to see that the  solutions of \eqref{b4a} and \eqref{b5a} are  critical points of the  functional
   \ber
    I_{\alpha\beta}(u, v)&=&\frac14\left(\frac1\alpha+\frac1\beta\right)\left(\|\nabla u\|_2^2+\|\nabla v\|_2^2\right)
    +\frac12\left(\frac1\alpha-\frac1\beta\right)\ito\nabla u\cdot\nabla v\ud x \nm\\
     &&+\alpha\ito\left(\re^{u_0+u}+\re^v-2\right)^2\ud x+\beta\ito\left(\re^{u_0+u}-\re^v\right)^2\ud x  \nm\\
       &&+\frac{4\pi n}{|\Omega|}\left(\frac1\alpha+\frac1\beta\right)\ito u \ud x
       +\frac{4\pi n}{|\Omega|}\left(\frac1\alpha-\frac1\beta\right)\ito v \ud x.\label{b6}
   \eer

In the following subsections we will  apply a constrained minimization  approach to find a first critical point  and the  mountain pass theorem to find a second critical point of the above functional, respectively.

\subsection{Constrained minimization}
Let $(u, v)$ be a solution of \eqref{b4a} and \eqref{b5a}, which   is also  a  solution of \eqref{b4} and \eqref{b5}.  Then  integrating these equations over  $\Omega$, we
obtain the following constraints
  \ber
  \alpha\ito\left(\re^{u_0+u}+\re^{v}-2\right)\left(\re^{u_0+u}+\re^{v}\right)\ud x
   +\beta\ito\left(\re^{u_0+u}-\re^{v}\right)^2\ud x+\frac{4\pi n}{\alpha}=0, \label{b6a}\\
    \alpha\ito\left(\re^{u_0+u}+\re^{v}-2\right)\left(\re^{u_0+u}-\re^{v}\right)\ud x
   +\beta\ito\left(\re^{2u_0+2u}-\re^{2v}\right)\ud x+\frac{4\pi n}{\beta}=0,\label{b6b}
  \eer
  or equivalently,
  \ber
  \ito\left(\re^{u_0+u}-1\right)\re^{u_0+u}\ud x-\gamma\ito\left(\re^{v}-1\right)\re^{u_0+v}\ud x+\frac{2\pi n}{\alpha\beta}=0,\label{b7}\\
  \ito\left(\re^v-1\right)\re^v\ud x-\gamma\ito\left(\re^{u_0+u}-1\right)\re^v\ud x+\frac{2\gamma\pi n}{\alpha\beta}=0,\label{b8}
  \eer
 where we define
  \be
  \gamma\equiv\frac{\beta-\alpha}{\beta+\alpha} \label{b9}
  \ee
   throughout the rest of the work.  Under our assumption on $\alpha, \beta$, that is, $\beta>\alpha>0$, we see that $0<\gamma<1$.

It can be checked that the constraints (\ref{b6a}) and (\ref{b6b}) are the quantized integrals (\ref{xq1}) and (\ref{xq2}) stated in Theorem \ref{thb1}.

From \eqref{b6a}, we see that
 \ber
 && \alpha\ito\left(\re^{u_0+u}+\re^v-2\right)^2\ud x+\beta\ito\left(\re^{u_0+u}-\re^v\right)^2\ud x
  \nm\\
  &&=2\alpha\left(\ito[1-\re^{u_0+u}]\ud x+\ito[1-\re^{v}] \ud x\right)-\frac{4\pi n}{\alpha}.\label{b8'}
 \eer

 We know that $W^{1,2}(\Omega)$ can be decomposed  as follows,
  \[
  W^{1,2}(\Omega)=\mathbb{R}\oplus  \dot{W}^{1,2}(\Omega),
  \]
  where
\[\dot{W}^{1,2}(\Omega)=\left\{w\in W^{1,2}(\Omega)\Bigg| \ito w\ud x=0\right\}\]
  is a closed subspace of $W^{1,2}(\Omega)$.

Then, we can decompose $u, v$ into the form
 \[u=u'+c_1,\quad v=v'+c_2, \]
where
  \[\ito u'\ud x=0,\quad \ito v'\ud x=0,\quad c_1=\frac{1}{|\Omega|}\ito u \ud x,   \quad c_2=\frac{1}{|\Omega|}\ito v\ud x. \]

Then \eqref{b7} and \eqref{b8} can be rewritten in  the  form
  \ber
   \re^{2c_1}\ito\re^{2u_0+2u'}\ud x-Q_1(u', v',  \re^{c_2})\re^{c_1}+\frac{2\pi n}{\alpha\beta}=0,\label{b10}\\
    \re^{2c_2}\ito\re^{2v'}\ud x-Q_2(u', v', \re^{c_1})\re^{c_2}+\frac{2\gamma\pi n}{\alpha\beta}=0,\label{b11}
  \eer
 where
  \ber
   Q_1(u', v',  \re^{c_2})&\equiv&(1-\gamma)\ito\re^{u_0+u'}\ud x+\gamma\re^{c_2}\ito\re^{u_0+u'+v'}\ud x,\label{b12}\\
    Q_2(u', v', \re^{c_1})&\equiv&(1-\gamma)\ito\re^{v'}\ud x+\gamma\re^{c_1}\ito\re^{u_0+u'+v'}\ud  x.\label{b13}
  \eer

Hence  the  equations \eqref{b10}--\eqref{b11} are solvable with respect to $c_1 $ and $c_2$ if and only if
 \ber
 (Q_1(u', v',  \re^{c_2}))^2&\ge&\frac{8\pi n}{\alpha\beta}\ito\re^{2u_0+2u'}\ud x,\label{b14}\\
  (Q_2(u', v', \re^{c_1}))^2&\ge&\frac{8\gamma\pi n}{\alpha\beta}\ito\re^{2v'}\ud x.\label{b15}
 \eer

From Proposition \ref{pr1}, we see that, for a solution $(u, v )$ of \eqref{b4a}--\eqref{b5a},   $u_0+u<0, v<0$, namely, $u_0+u'+c_1<0, v'+c_2<0$. Then from \eqref{b14} we obtain
 \[
  \frac{8\pi n}{\alpha\beta}\ito\re^{2u_0+2u'}\ud x\le\left(\ito\re^{u_0+u'}\ud x\right)^2\le |\Omega|\ito\re^{2u_0+2u'}\ud x,
 \]
which gives a necessary condition for the existence of solutions to \eqref{b4a}--\eqref{b5a}
 \be
  \alpha\beta\ge\frac{8\pi n}{|\Omega|}. \label{b16}
 \ee
Then we get the first conclusion of Theorem \ref{thb1}.

Now we take  the following constraints
 \ber
  \left(\ito\re^{u_0+u'}\ud x\right)^2&\ge&\frac{8\pi  n}{(1-\gamma)^2\alpha\beta}\ito\re^{2u_0+2u'}\ud x,\label{b17}\\
    \left(\ito\re^{v'}\ud x\right)^2&\ge&\frac{8\gamma\pi  n}{(1-\gamma)^2\alpha\beta}\ito\re^{2v'}\ud x.\label{b18}
 \eer
 We  introduce the following   admissible set
  \be
  \mathcal{A}=\left\{(u', v')\in \dot{W}^{1,2}(\Omega)\times\dot{W}^{1,2}(\Omega)\Big|(u', v') \quad \text{satisfies} \quad\eqref{b17}-\eqref{b18}\right\}. \label{b19}
  \ee

Thus,  for any  $(u', v')\in\mathcal{A}$,  we can find  a solution of the equations \eqref{b10}--\eqref{b11} with respect to $c_1$ and $c_2$  by solving the following equations
 \ber
 \re^{c_1}&=&\frac{ Q_1(u', v',  \re^{c_2})+\sqrt{[Q_1(u', v',  \re^{c_2})]^2-\frac{8\pi  n}{\alpha\beta}\ito\re^{2u_0+2u'}\ud x}}{2\ito\re^{2u_0+2u'}\ud x}
 \nm\\&\equiv& g_1(\re^{c_2}), \label{b20}\\
  \re^{c_2}&=&\frac{ Q_2(u', v',  \re^{c_1})+\sqrt{[Q_2(u', v',  \re^{c_1})]^2-\frac{8\gamma\pi  n}{\alpha\beta}\ito\re^{2v'}\ud x}}{2\ito\re^{2v'}\ud x}
 \nm\\&\equiv& g_2(\re^{c_1}).\label{b21}
 \eer

Indeed,  letting
 \[ F(X)\equiv X-g_1(g_2(X)),  \]
 we can solve \eqref{b20}--\eqref{b21}  by finding the zeros of the function $F(\cdot)$. Therefore, it is sufficient to prove the following proposition.
 \begin{proposition}\label{pr2}
 For any $(u', v')\in   \mathcal{A}$, the equation
  \[F(X)=X-g_1(g_2(X))=0\]
   admits a unique positive solution $X_0$.
 \end{proposition}

By this proposition, for any  $(u', v')\in\mathcal{A}$,  we can get a solution of  \eqref{b10} and \eqref{b11} with respect to $ c_1, c_2$.

 {\bf Proof of the Proposition \ref{pr2}. }   \quad By  \eqref{b20} and \eqref{b21} it is easy to see that
 \be
 g_i(X)>0, \quad \forall\, X\ge0, \, i=1,2.\label{b28}
 \ee
 Then, we see that $F(0)=-g_1(g_2(0))<0$. We check that
 \ber
  \frac{\ud g_1(X)}{\ud X}&=& \frac{\gamma g_1(X)\ito\re^{u_0+u'+v'}\ud x}{\sqrt{[Q_1(u', v', X)]^2-\frac{8\pi n}{\alpha\beta}\ito\re^{2u_0+2u'}\ud x}}, \label{b26}\\
  \frac{\ud g_2(X)}{\ud X}&=& \frac{\gamma g_2(X)\ito\re^{u_0+u'+v'}\ud x}{\sqrt{[Q_2(u', v', X)]^2-\frac{8\gamma\pi n}{\alpha\beta}\ito\re^{2v'}\ud x}},\label{b27}
 \eer
  which are all positive. Then, we see that $g_i(X)$ $(i=1, 2)$ is  strictly increasing for all $X>0$.

  After a direct computation, we obtain
    \berr
    \lim\limits_{X\to+\infty}\frac{g_1(X)}{X}=\frac{\gamma\ito\re^{u_0+u'+v'}\ud x}{\ito\re^{2u_0+2u'}\ud x}, \\
    \lim\limits_{X\to+\infty}\frac{g_2(X)}{X}=\frac{\gamma\ito\re^{u_0+u'+v'}\ud x}{\ito\re^{2v'}\ud x},
    \eerr
 from which it follows that
   \berr
  \lim\limits_{X\to+\infty}\frac{F(X)}{X}&=&1-\frac{\gamma^2\left(\ito\re^{u_0+u'+v'}\ud x\right)^2}{\ito\re^{2u_0+2u'}\ud x\ito\re^{2v'}\ud x}
  \ge1-\gamma^2>0.
   \eerr
Therefore, we  have \[\lim\limits_{X\to+\infty}F(X)=+\infty.\]

Noting that $F(0)<0$, then we conclude that the  equation $F(X)=0$ has at least one  solution $X_0>0$.

Next we show that the solution is also unique.
From \eqref{b26} and \eqref{b27} and \eqref{b17} and \eqref{b18} we obtain
 \berr
  \frac{\ud F(X)}{\ud X}&=&1-\frac{g_1(g_2(X))g_2(X)\gamma^2\left(\ito\re^{u_0+u'+v'}\ud x\right)^2}
  {\sqrt{[Q_1(u', v', g_2(X))]^2-\frac{8\pi n}{\alpha\beta}\ito\re^{2u_0+2u'}\ud x}\sqrt{[Q_2(u', v', X)]^2-\frac{8\gamma\pi n}{\alpha\beta}\ito\re^{2v'}\ud x}}
\\
  &>&1-\frac{g_1(g_2(X))}{X}=\frac{F(X)}{X}.
 \eerr
 Then we have  \[  \frac{\ud }{\ud X}\left(\frac{F(X)}{X}\right)>0, \]
  which says that $\frac{F(X)}{X}$ is strictly increasing  for $X>0$. As a result, $F(X)$ is strictly increasing for $X>0$. Then $F(X)$ has a unique zero point.
 Then the proof of Proposition \ref{pr2} is complete.

 By the above discussion we see that, for any $(u', v')\in \mathcal{A}$,  we can get pair $(c_1(u',  v'), c_2(u',  v'))$ given by
 \eqref{b20}--\eqref{b21},  which solves \eqref{b10}--\eqref{b11}, such that  $(u, v)$ defined by
 \[ u=u'+c_1(u',  v'), \quad  v=v'+c_2(u',  v') \]
 satisfies  \eqref{b6a}--\eqref{b6b}.

In what follows we consider the  minimization problem
  \ber
   \min \left\{J_{\alpha\beta}(u', v')\big|\, (u', v')\in \mathcal{A}\right\},\label{5.30}
  \eer
  where  $J_{\alpha\beta}(u', v')$ is defined by
 \berr
 J_{\alpha\beta}(u', v')=I_{\alpha\beta}(u'+c_1(u',  v'), v'+c_2(u', v')),
 \eerr
 $(c_1(u',  v'), c_2(u',  v'))$ is  given by \eqref{b20}--\eqref{b21}.  From \eqref{b6} and \eqref{b8'}, we see that
 \ber
   J_{\alpha\beta}(u', v')&=&\frac14\left(\frac1\alpha+\frac1\beta\right)\left(\|\nabla u'\|_2^2+\|\nabla v'\|_2^2\right)
    +\frac12\left(\frac1\alpha-\frac1\beta\right)\ito\nabla u'\cdot\nabla v'\ud x \nm\\
    &&+2\alpha\left(\ito\big[1-\re^{u_0+u'}\re^{c_1}\big]\ud x+\ito\big[1-\re^{v'}\re^{c_2}\big]\ud x\right)-\frac{4\pi n}{\alpha} \nm\\
    &&+4\pi n\left(\frac1\alpha+\frac1\beta\right)c_1+4\pi n\left(\frac1\alpha-\frac1\beta\right)c_2.\label{b29}
 \eer

It is easy to check that the functional  $J_{\alpha\beta}$ is Frech\'{e}t differentiable in the interior of $\mathcal{A}$.
Moreover, if $(u' ,v')$ is an interior  critical point of $J_{\alpha\beta}$ in $\mathcal{A}$,  then $ (u'+c_1(u',  v'), v'+c_2(u', v'))$
gives a critical point for $I_{\alpha\beta}$.

In  what follows, we will prove that the functional  $J_{\alpha\beta}$ is bounded from  below and the problem \eqref{5.30} admits interior minimum.

\begin{lemma}\label{lem1}
 For any $(u', v')\in \mathcal{A}$, we have
  \be \re^{c_1}\ito \re^{u_0+u'}\ud x\le |\Omega|, \quad \re^{c_2}\ito \re^{v'}\ud x\le |\Omega|.\label{b30}
  \ee
\end{lemma}
\begin{remark}\label{rmk1}
 It follows from the  Jensen inequality and \eqref{b30} that
  \[\re^{c_1}\le 1,\quad \re^{c_2}\le 1.\]
\end{remark}

{\bf Proof.}  \quad From  \eqref{b20} and \eqref{b21}, we  obtain
 \ber
  \re^{c_1}\le \frac{Q_1(u', v',  \re^{c_2})}{\ito \re^{2u_0+2u'}\ud x},\label{b31}\\
    \re^{c_2}\le \frac{Q_2(u', v',  \re^{c_1})}{\ito \re^{ 2v'}\ud  x}. \label{b32}
 \eer
Then it follows from \eqref{b31}, \eqref{b32}, \eqref{b12}, \eqref{b13},  and the  H\"{o}lder inequality that
 \berr
  \re^{c_1}&\le& \frac{(1-\gamma)\ito\re^{u_0+u'}\ud x}{\ito \re^{2u_0+2u'}\ud x}
   +\frac{\gamma(1-\gamma)\ito\re^{u_0+u'+v'}\ud x\ito\re^{v'}\ud x}{\ito \re^{2u_0+2u'}\ud x\ito \re^{2v'}\ud x}
    +\frac{\gamma^2\left(\ito\re^{u_0+u'+v'}\ud x\right)^2}{\ito \re^{2u_0+2u'}\ud x\ito \re^{2v'}\ud x}\re^{c_1}\nm\\
     &\le& \frac{(1-\gamma)\ito\re^{u_0+u'}\ud x}{\ito \re^{2u_0+2u'}\ud x}
    +\frac{\gamma(1-\gamma)\ito\re^{u_0+u'+v'}\ud x\ito\re^{v'}\ud x}{\ito \re^{2u_0+2u'}\ud x\ito \re^{2v'}\ud x}
    +\gamma^2\re^{c_1},
 \eerr
which enables us to conclude that
 \ber
 \re^{c_1}\le \frac{1}{1+\gamma}\left(\frac{\ito\re^{u_0+u'}\ud x}{\ito \re^{2u_0+2u'}\ud x}
    +\frac{\gamma\ito\re^{u_0+u'+v'}\ud x\ito\re^{v'}\ud x}{\ito \re^{2u_0+2u'}\ud x\ito \re^{2v'}\ud  x}\right). \label{b33}
 \eer
 Similarly, we have
 \ber
 \re^{c_2}\le \frac{1}{1+\gamma}\left(\frac{\ito\re^{v'}\ud x}{\ito \re^{2v'}\ud x}
    +\frac{\gamma\ito\re^{u_0+u'+v'}\ud x\ito\re^{u_0+u'}\ud x}{\ito \re^{2u_0+2u'}\ud x\ito \re^{2v'}\ud x}\right). \label{b34}
 \eer
Using \eqref{b33} and \eqref{b34} and  the H\"{o}lder inequality, we
have
 \berr
 \re^{c_1}\ito\re^{u_0+u'}\ud x&\le& \frac{1}{1+\gamma}\left(\frac{\left[\ito\re^{u_0+u'}\ud x\right]^2}{\ito \re^{2u_0+2u'}\ud x}
    +\frac{\gamma\ito\re^{u_0+u'+v'}\ud x\ito\re^{u_0+u'}\ud x\ito\re^{v'}\ud x}{\ito \re^{2u_0+2u'}\ud x\ito \re^{2v'}\ud x}\right)\\
     &\le&|\Omega|, \\
     \re^{c_2}\ito\re^{v'}\ud x&\le& \frac{1}{1+\gamma}\left(\frac{\left[\ito\re^{v'}\ud x\right]^2}{\ito \re^{2v'}\ud x}
    +\frac{\gamma\ito\re^{u_0+u'+v'}\ud x\ito\re^{u_0+u'}\ud x\ito\re^{v'}\ud x}{\ito \re^{2u_0+2u'}\ud x\ito \re^{2v'}\ud x}\right)\\
     &\le&|\Omega|.
 \eerr
Thus the lemma follows.

Estimates of the type contained in the following lemma were observed first in \cite{nota}.
 \begin{lemma}\label{lem2}
  For any $(u', v')\in \mathcal{A}$ and $s\in (0, 1)$,  it holds
    \ber
   \ito\re^{u_0+u'}\ud x&\le& \left(\frac{[1-\gamma]^2\alpha\beta}{8\pi  n}\right)^{\frac{1-s}{s}}
   \left(\ito\re^{su_0+su'}\ud x\right)^{\frac1s}, \label{b35}\\
     \ito\re^{v'}\ud x&\le& \left(\frac{[1-\gamma]^2\alpha\beta}{8\gamma\pi  n}\right)^{\frac{1-s}{s}}
   \left(\ito\re^{sv'}\ud x\right)^{\frac1s}. \label{b36}
    \eer

 \end{lemma}
{\bf Proof.}  \quad  Let $s\in (0, 1)$, $a=\frac{1}{2-s}$ such that
 $sa+2(1-a)=1$. Then  using the H\"{o}lder inequality  and \eqref{b17} we have
   \berr
   \ito\re^{u_0+u'}\ud x&\le&
   \left(\ito\re^{su_0+su'}\ud x\right)^a\left(\ito\re^{2u_0+2u'}\ud x\right)^{1-a}\\
   &\le& \left(\frac{[1-\gamma]^2\alpha\beta}{8\pi n}\right)^{1-a}\left(\ito\re^{su_0+su'}\ud x\right)^a\left(\ito\re^{u_0+u'}\ud x\right)^{2(1-a)},
   \eerr
which implies
 \berr
   \ito\re^{u_0+u'}\ud x
   &\le& \left(\frac{[1-\gamma]^2\alpha\beta}{8\pi n}\right)^{\frac{1-a}{2a-1}}\left(\ito\re^{su_0+su'}\ud x\right)^{\frac{a}{2a-1}}\\
 &=&\left(\frac{[1-\gamma]^2\alpha\beta}{8\pi  n}\right)^{\frac{1-s}{s}} \left(\ito\re^{su_0+su'}\ud x\right)^{\frac1s}.
   \eerr
   Analogously,  we can obtain \eqref{b36}.

Next we show that the functional $J_{\alpha\beta}$ is coercive and bounded from below  on $ \mathcal{A}$.
  To this end, we will use the  Trudinger--Moser inequality (see \cite{font})
 \be
 \ito \re^{w}\ud x \le C_1\exp\left(\frac{1}{16\pi}\|\nabla w\|_2^2\right), \quad \forall\, w\in\dot{W}^{1,2}(\Omega),\label{b37}
 \ee
 where $C_1$ is a positive constant.

 \begin{lemma}\label{lem3}
   For any $(u', v')\in \mathcal{A}$, the functional
   $J_{\alpha\beta}$ satisfies
    \ber
      J_{\alpha\beta}(u', v') \ge \frac{1}{4\beta}\left(\|\nabla u'\|_2^2+\|\nabla v'\|_2^2\right)
      -C_{\alpha\beta}, \label{b38}
    \eer
    where
    \ber
     C_{\alpha\beta}&\equiv& \frac{8\pi n^2(\alpha+\beta)}{\alpha^2}\left(\ln\alpha\beta+\ln C_1+\ln\frac{[1-\gamma]^2}{8\pi n}\right)
     -\frac{8\pi n}{\alpha}\ln\frac{1-\gamma}{2} \nm\\
     &&+4\pi n\left(\frac1\alpha+\frac1\beta\right)\max\limits_{x\in\Omega}u_0
       -4\pi n^2\left(\frac1\alpha-\frac1\beta\right)\left(1+\frac\beta\alpha\right)\ln\gamma.\label{b39}
    \eer
 \end{lemma}
 {\bf Proof.}\quad From \eqref{b20} and \eqref{b21}, we see that
  \[
  \re^{c_1}\ge \frac{(1-\gamma)\ito\re^{u_0+u'}\ud x}{2\ito \re^{2u_0+2u'}\ud x},
   \quad \re^{c_2}\ge \frac{(1-\gamma)\ito\re^{v'}\ud x}{2\ito \re^{2v'}\ud x}.
  \]
Then by \eqref{b17} and \eqref{b18}, we obtain
 \berr
 \re^{c_1}\ge\frac{4\pi n}{(1-\gamma)\alpha\beta\ito\re^{u_0+u'}\ud x}, \quad
   \re^{c_2}\ge\frac{4\gamma\pi n}{(1-\gamma)\alpha\beta\ito\re^{v'}\ud x},
 \eerr
which lead to
 \ber
  c_1&\ge& \ln\frac{4\pi n}{1-\gamma}-\ln\alpha\beta-\ln\ito\re^{u_0+u'}\ud x,\label{b40}\\
  c_2&\ge& \ln\frac{4\gamma\pi n}{1-\gamma}-\ln\alpha\beta-\ln\ito\re^{v'}\ud x.\label{b41}
 \eer

For any $s\in(0, 1)$,  using  Lemma \ref{lem2} and the Trudinger--Moser inequality \eqref{b37}, we have
 \ber
  \ln\ito\re^{u_0+u'}\ud x&\le& \frac{1-s}{s}\left(\ln\frac{[1-\gamma]^2}{8\pi n}+\ln\alpha\beta\right)+\frac1s\ln\ito\re^{su_0+su'}\ud x\nm\\
   &\le &\frac{s}{16\pi}\|\nabla u'\|_2^2+\frac{1-s}{s}\left(\ln\frac{[1-\gamma]^2}{8\pi n}+\ln\alpha\beta\right)+\max\limits_{x\in\Omega}u_0+\frac1s\ln
   C_1,\label{b42}
   \\
     \ln\ito\re^{v'}\ud x&\le& \frac{1-s}{s}\left(\ln\frac{[1-\gamma]^2}{8\gamma\pi n}+\ln\alpha\beta\right)+\frac1s\ln\ito\re^{sv'}\ud x\nm\\
   &\le &\frac{s}{16\pi}\|\nabla v'\|_2^2+\frac{1-s}{s}\left(\ln\frac{[1-\gamma]^2}{8\gamma\pi n}+\ln\alpha\beta\right)+\frac1s\ln C_1.\label{b43}
 \eer
Then, from \eqref{b29}, \eqref{b8'}, \eqref{b40}--\eqref{b43}, we see that
  \ber
   J_{\alpha\beta}(u', v')&\ge& \left(\frac{1}{2\beta}-\frac{s n}{4}\left[\frac1\alpha+\frac1\beta\right]\right)\|\nabla u'\|_2^2
     + \left(\frac{1}{2\beta}-\frac{sn}{4}\left[\frac1\alpha-\frac1\beta\right]\right)\|\nabla v'\|_2^2\nm\\
     &&+\frac{8\pi n}{\alpha}\ln\frac{1-\gamma}{2}-4\pi n\left(\frac1\alpha+\frac1\beta\right)\max\limits_{x\in\Omega}u_0\nm\\
      &&-\frac{8\pi n}{s\alpha}\left(\ln\alpha\beta+\ln C_1+\ln\frac{[1-\gamma]^2}{8\pi n}\right)
      +\frac{4\pi n}{s}\left(\frac1\alpha-\frac1\beta\right)\ln\gamma.\label{b44}
  \eer
 Now by taking    \[s=\frac{\alpha}{n(\alpha+\beta)}, \]  in \eqref{b44}, we get \eqref{b38}.

It is easy to see that $ J_{\alpha\beta}(u', v')$ is weakly lower semi-continuous on $\mathcal{A}$. Then  by lemma \ref{lem3} we infer
that the infimum of $J_{\alpha\beta}(u', v')$ can be attained on $\mathcal{A}$.

In  the sequel we will show  that, when $\alpha, \beta$ satisfy $\beta>\alpha>0$, \eqref{b03}, and  $\alpha$ is sufficiently large,   a minimizer can
only be an interior point of $\mathcal{A}$.

 \begin{lemma}\label{lem4}
  The functional $J_{\alpha\beta}$ satisfies
   \be
   \inf\limits_{(u', v')\in \partial\mathcal{A}}J_{\alpha\beta}(u', v')
    \ge 2|\Omega|\alpha-\frac{16\pi n}{(1+\gamma)(1-\gamma)^2\alpha} -\frac{4\gamma\sqrt{2\pi n|\Omega|}}{1-\gamma^2}-C_{\alpha\beta},\label{b45}
   \ee
   where $C_{\alpha\beta}$ is defined by \eqref{b39}.
 \end{lemma}

{\bf Proof.} \quad  On the boundary of $\mathcal{A}$, we have
  \ber
  \left(\ito\re^{u_0+u'}\ud x\right)^2&=&\frac{8\pi  n}{(1-\gamma)^2\alpha\beta}\ito\re^{2u_0+2u'}\ud x\label{b46}
  \eer
   or
   \ber \left(\ito\re^{v'}\ud x\right)^2&=&\frac{8\gamma\pi  n}{(1-\gamma)^2\alpha\beta}\ito\re^{2v'}\ud x.\label{b47}
 \eer

If \eqref{b46} holds, using \eqref{b33} and the H\"{o}lder inequality,  we obtain
 \berr
 \re^{c_1}\ito\re^{u_0+u'}\ud x&\le& \frac{1}{1+\gamma}\left(\frac{\left[\ito\re^{u_0+u'}\ud x\right]^2}{\ito \re^{2u_0+2u'}\ud x}
    +\frac{\gamma\ito\re^{u_0+u'+v'}\ud x\ito\re^{u_0+u'}\ud x\ito\re^{v'}\ud x}{\ito \re^{2u_0+2u'}\ud x\ito \re^{2v'}\ud x}\right)\\
     &\le&\frac{8\pi n}{(1+\gamma)(1-\gamma)^2\alpha\beta}+\frac{2\gamma\sqrt{2\pi n|\Omega|}}{(1-\gamma^2)\sqrt{\alpha\beta}}\\
     &\le&\frac{8\pi n}{(1+\gamma)(1-\gamma)^2\alpha^2}+\frac{2\gamma\sqrt{2\pi n|\Omega|}}{(1-\gamma^2)\alpha},
  \eerr
  which leads to
   \berr
   2\alpha\left(\ito\big[1-\re^{u_0+u'}\re^{c_1}\big]\ud x+\ito\big[1-\re^{v'}\re^{c_2}\big]\ud x\right)
   \ge2|\Omega|\alpha-\frac{16\pi n}{(1+\gamma)(1-\gamma)^2\alpha} -\frac{4\gamma\sqrt{2\pi n|\Omega|}}{1-\gamma^2}.
   \eerr

Therefore, using  similar estimates for $c_1, c_2$ as in Lemma \ref{lem3}, on the boundary of $\mathcal{A}$  we have
  \berr
  J_{\alpha\beta}(u', v')\ge2|\Omega|\alpha-\frac{16\pi n}{(1+\gamma)(1-\gamma)^2\alpha} -\frac{4\gamma\sqrt{2\pi n|\Omega|}}{1-\gamma^2}-C_{\alpha\beta},
  \eerr
  where $C_{\alpha\beta}$ is defined by \eqref{b39}.  Then the proof of Lemma \ref{lem4} is complete.

As a test function, we  use,  as in \cite{nota},  the solution characterized by   Tarantello \cite{taran96}. Namely, from   \cite{taran96},    we know  that, for  $\lambda>0$ sufficiently large, there exists a solution $w_{\lambda}$ of the equation
   \ber
   \Delta w=\lambda\re^{u_0+w}\big(\re^{u_0+w}-1\big)+\frac{8\pi n}{|\Omega|}, \label{b48}
   \eer
satisfying   $w_\lambda=c_\lambda+w_\lambda'$, $c_\lambda=\frac{1}{|\Omega|}\ito w_\lambda\ud  x$, $\ito w'_\lambda\ud x=0$, such that $u_0+w_\lambda<0$ in
$\Omega$,  $c_\lambda\to 0$, and $w'_\lambda\to -u_0$ pointwise, as $\lambda\to +\infty$.

In view of $\re^{u_0}\in L^{\infty}(\Omega)$ and the dominated  convergence theorem, we have
 \[\re^{u_0+w_\lambda'}\to1 \quad \text{strongly in }\quad L^p(\Omega) \quad \text{for any } p\ge 1, \]
as $\lambda\to +\infty$. In particular,
  \[\ito\re^{2u_0+2w_\lambda'}\ud x\to |\Omega|,  \]
 as $\lambda\to +\infty$. Therefore, for $\alpha_0$ large and fixed
 $\vep\in (0, 1)$, we can find $\lambda_\vep$ to ensure that $(w'_{\lambda_\vep}, 0)\in
 \mathcal{A}$ for every $\alpha>\alpha_0$ and
   \be
   \frac{(1-\gamma^2)|\Omega|}{\ito\re^{2u_0+2w'_{\lambda_\vep}}\ud x-\gamma^2|\Omega|}\ge1-\vep. \label{b49}
   \ee

 By the Jensen inequality,
   \[ \ito \re^{u_0+w'_{\lambda_\vep}}\ud x\ge |\Omega|.\]

 Then in view of \eqref{b20} and \eqref{b21} we   get
  \ber
   \re^{c_1(w'_{\lambda_\vep}, 0)}&\ge& \frac{Q_1(w'_{\lambda_\vep}, 0, \re^{c_2(w'_{\lambda_\vep}, 0)})}{2\ito\re^{2u_0+2w'_{\lambda_\vep}}\ud x}\left(1+\sqrt{1-\frac{8\pi n\ito\re^{2u_0+2w'_{\lambda_\vep}}\ud x}{\alpha\beta Q^2_1(w'_{\lambda_\vep}, 0, \re^{c_2(w'_{\lambda_\vep}, 0)})}}\right)\nm\\
   &\ge&  \frac{Q_1(w'_{\lambda_\vep}, 0, \re^{c_2(w'_{\lambda_\vep}, 0)})}{\ito\re^{2u_0+2w'_{\lambda_\vep}}\ud x}-\frac{4\pi n}{\alpha\beta Q_1(w'_{\lambda_\vep}, 0, \re^{c_2(w'_{\lambda_\vep}, 0)})} \nm\\
   &\ge& \frac{(1-\gamma+\gamma\re^{c_2(w'_{\lambda_\vep}, 0)})|\Omega|}{\ito\re^{2u_0+2w'_{\lambda_\vep}}\ud x}-\frac{4\pi n}{\alpha\beta(1-\gamma)|\Omega|}. \label{b51}
  \eer
 Similarly, we have
  \ber
   \re^{c_2(w'_{\lambda_\vep}, 0)}&\ge&  1-\gamma +\gamma\re^{c_1(w'_{\lambda_\vep}, 0)}-\frac{4\gamma\pi n}{\alpha\beta(1-\gamma)|\Omega|}. \label{b52}
  \eer

Therefore  it follows from  \eqref{b51} and \eqref{b52} that
  \berr
   \re^{c_1(w'_{\lambda_\vep}, 0)}&\ge&
    \frac{\gamma|\Omega|}{\ito\re^{2u_0+2w'_{\lambda_\vep}}\ud x}\left(1-\gamma +\gamma\re^{c_1(w'_{\lambda_\vep}, 0)}-\frac{4\gamma\pi n}{\alpha\beta[1-\gamma]|\Omega|}\right)
     \nm\\ &&+\frac{(1-\gamma)|\Omega|}{\ito\re^{2u_0+2w'_{\lambda_\vep}}\ud x} -\frac{4\pi n}{\alpha\beta(1-\gamma)|\Omega|} \nm\\
     &\ge& \frac{(1-\gamma^2)|\Omega|}{\ito\re^{2u_0+2w'_{\lambda_\vep}}\ud x}
      +\frac{ \gamma^2|\Omega| \re^{c_1(w'_{\lambda_\vep}, 0)}}{\ito\re^{2u_0+2w'_{\lambda_\vep}}\ud x}- \frac{4\pi n(1+\gamma^2)}{\alpha\beta(1-\gamma)|\Omega|},
  \eerr
which implies
  \ber
   \re^{c_1(w'_{\lambda_\vep}, 0)}&\ge&  \frac{(1-\gamma^2)|\Omega|}{\ito\re^{2u_0+2w'_{\lambda_\vep}}\ud x-\gamma^2|\Omega|}
  - \frac{ 4\pi n(1+\gamma^2) \ito\re^{2u_0+2w'_{\lambda_\vep}}\ud x}{\alpha\beta(1-\gamma)|\Omega| \left(\ito\re^{2u_0+2w'_{\lambda_\vep}}\ud x-\gamma^2|\Omega|\right)}
  \nm\\ &\ge& \frac{(1-\gamma^2)|\Omega|}{\ito\re^{2u_0+2w'_{\lambda_\vep}}\ud x-\gamma^2|\Omega|}
    -\frac{4\pi n(1+\gamma^2)}{\alpha\beta(1+\gamma)(1-\gamma)^2|\Omega|}.\label{b53}
 \eer

Similarly, we get

  \ber
   \re^{c_2(w'_{\lambda_\vep}, 0)} &\ge& \frac{(1-\gamma^2)|\Omega|}{\ito\re^{2u_0+2w'_{\lambda_\vep}}\ud x-\gamma^2|\Omega|}
    -\frac{8\pi n\gamma}{\alpha\beta(1+\gamma)(1-\gamma)^2|\Omega|}.\label{b54}
 \eer

Then,  combining \eqref{b53}, \eqref{b54}, and \eqref{b49},  we conclude that for all $\alpha>\alpha_0$,
 \berr
 \re^{c_1(w'_{\lambda_\vep}, 0)} \ge 1-\vep-\frac{4\pi n(1+\gamma^2)}{\alpha\beta(1+\gamma)(1-\gamma)^2|\Omega|},\\
  \re^{c_2(w'_{\lambda_\vep}, 0)} \ge 1-\vep-\frac{8\pi n\gamma}{\alpha\beta(1+\gamma)(1-\gamma)^2|\Omega|}.
 \eerr

Therefore, we obtain that for all $\alpha>\alpha_0$,
  \ber
   \ito\left(1-\re^{c_1(w'_{\lambda_\vep}, 0)}\re^{u_0+w'_{\lambda_\vep}}\right)\ud x&\le& \vep |\Omega|+\frac{4\pi n(1+\gamma^2)}{\alpha\beta(1+\gamma)(1-\gamma)^2},\label{b55} \\
   \ito\left(1-\re^{c_2(w'_{\lambda_\vep}, 0)}\right)\ud x&\le& \vep |\Omega|+\frac{8\pi n\gamma}{\alpha\beta(1+\gamma)(1-\gamma)^2}. \label{b56}
  \eer

\begin{lemma}\label{lem5}
  Assume $\beta>\alpha>0$ and that \eqref{b03} holds.  There exists a positive constant
  $M_\sigma$ such that,  when $\alpha>M_\sigma$, we have
    \be
    J_{\alpha\beta}(w'_{\lambda_\vep}, 0)-  \inf\limits_{(u', v')\in \partial\mathcal{A}}J_{\alpha\beta}(u', v')
    <-1. \label{b57}
    \ee
\end{lemma}

 {\bf Proof.} \quad Fix $\vep\in(0, \frac14)$, we choose $w'_{\lambda_\vep}$ that 
 satisfies \eqref{b55}--\eqref{b56}.  Noting that $\re^{c_1}\le 1,\, \re^{c_2}\le 1$, then by \eqref{b29}, \eqref{b55} and \eqref{b56}  we have
   \ber
   J_{\alpha\beta}(w'_{\lambda_\vep}, 0)&\le& \frac14\left(\frac1\alpha+\frac1\beta\right)\|\nabla w'_{\lambda_\vep}\|_2^2+4\vep\alpha|\Omega|
   +\frac{8\pi n(1+\gamma)}{(1-\gamma)^2\beta}-\frac{4\pi  n}{\alpha}\nm\\
   &\le&4\vep\alpha|\Omega|+\frac{C_\vep}{4}\left(\frac1\alpha+\frac1\beta\right)
    +\frac{8\pi n(1+\gamma)}{(1-\gamma)^2\beta}-\frac{4\pi  n}{\alpha},\label{b58}
   \eer
where $C_\vep$ is a positive constant depending only on $\vep$.

Then by Lemma \ref{lem4} we have
 \ber
 &&J_{\alpha\beta}(w'_{\lambda_\vep}, 0)-\inf\limits_{(u', v')\in \partial\mathcal{A}}J_{\alpha\beta}(u', v')\nm\\
  &&\le 2\alpha|\Omega|(2\vep-1)+\frac{C_\vep}{4}\left(\frac1\alpha+\frac1\beta\right)+\frac{8\pi n(1+\gamma)}{(1-\gamma)^2\beta}
  \nm \\&&\quad +\frac{16\pi n}{(1+\gamma)(1-\gamma)^2\alpha}+\frac{4\gamma\sqrt{2\pi n|\Omega|}}{1-\gamma^2}+C_{\alpha\beta}, \label{b59}
 \eer
where $C_{\alpha\beta}$ is defined by \eqref{b39}.

 By the assumption $\beta>\alpha>0$ and \eqref{b03}, we  easily deduce the following estimate
  \ber
  \frac{2}{1+\sigma}\le 1-\gamma<1,\label{b59a}\\
  C_{\alpha\beta}\le C_\sigma\left(\frac{\ln\alpha}{\alpha}+1\right)\label{b59b}
  \eer
  with a suitable constant $C_\sigma>0$ depending on $\sigma>1$ only.  By using the above estimates in \eqref{b59}, we obtain the desired  conclusion. 
  
 Using Lemma \ref{lem3} and \ref{lem5}, we  infer that under the assumption $\beta>\alpha>0$ and \eqref{b03}, when
 $\alpha$ is sufficiently large  the  problem \eqref{5.30}  admits a  minimizer $(u'_\alpha, v'_\alpha)$
 , which lies  in the interior of $\mathcal{A}$.  Consequently,
  \[ DJ_{\alpha\beta}(u'_\alpha, v'_\alpha)=0,\]
   and
    \ber
     u_\alpha=u'_\alpha+c_1(u'_\alpha, v'_\alpha), \quad v_\alpha=v'_\alpha+c_2(u'_\alpha, v'_\alpha) \label{b60}
    \eer
gives rise to  a critical point for the functional $I_{\alpha\beta}$, namely,  a weak solution to \eqref{b4a} and \eqref{b5a}.

In what follows we study the behavior of the solution  given by \eqref{b60}.

 \begin{lemma} \label{lem6}
 Let $(u_\alpha, v_\alpha)$ be defined by \eqref{b60}. Then
   \ber
  \re^{u_0+u_\alpha}\to1, \quad \re^{v_\alpha}\to1 \label{b61}
   \eer
    as $\alpha\to +\infty$ pointwise  a.e.  in $\Omega$,  and in  $L^p(\Omega)$ for any $p\ge1$.
 \end{lemma}

{\bf Proof.} \quad Using \eqref{b58}, we obtain that for any $\vep>0$ (small) there exists $\alpha_\vep>0$ such that when $\alpha>\alpha_\vep$,
 \be
 \inf\limits_{(u', v')\in \mathcal{A}}J_{\alpha\beta}(u', v')\le 4\vep\alpha|\Omega|+\frac{C_\vep}{4}\left(\frac1\alpha+\frac1\beta\right)
    +\frac{8\pi n(1+\gamma)}{(1-\gamma)^2\beta}-\frac{4\pi  n}{\alpha}. \label{b63}
 \ee

By a similar argument as in Lemma \ref{lem3} we have
 \ber
   &&\inf\limits_{(u', v')\in \mathcal{A}}J_{\alpha\beta}(u', v')\nm\\
   &&=J_{\alpha\beta}(u'_\alpha, v'_\alpha) \nm\\
    &&\ge \alpha \ito\left(\re^{u_0+u_\alpha}+\re^{v_\alpha}-2\right)^2\ud x+\beta\ito\left(\re^{u_0+u_\alpha}-\re^{v_\alpha}\right)^2\ud x-C_{\alpha\beta}  \nm\\
    &&\ge 2\alpha \left(\ito\left[\re^{u_0+u_\alpha}-1\right]^2\ud x+\ito\left[\re^{v_\alpha}-1\right]^2\ud x\right) -C_{\alpha\beta}, \label{b64}
 \eer
 where $C_{\alpha\beta}$ is defined by \eqref{b39}.

Then  by means of \eqref{b59a}, \eqref{b59b},  \eqref{b63}, and \eqref{b64}, we see that
  \berr
    \limsup\limits_{\alpha\to+\infty}\left(\ito\left[\re^{u_0+u_\alpha}-1\right]^2\ud x
     +\ito\left[\re^{v_\alpha}-1\right]^2\ud x\right)\le 2|\Omega|\vep, \quad \forall \vep>0,
  \eerr
which enables us to conclude  that,
  \[\re^{u_0+u_\alpha}\to 1, \quad  \re^{v_\alpha}\to 1,  \]
 in $L^2(\Omega)$,  as $\alpha\to+\infty$.    Since by Proposition \ref{pr1},  we also  know that
  \[\re^{u_0+u_\alpha}<1,\,  \re^{v_\alpha}<1\] in $\Omega$. Then we have
 \[ \re^{u_0+u_\alpha}\to1, \quad \re^{v_\alpha}\to1\]
 pointwise a.e.  in $\Omega$ as $\alpha\to+\infty$. At this point, we may complete the proof of Lemma \ref{lem6}  by  dominated convergence theorem.

 In order to get a secondary solution of \eqref{b4a} and \eqref{b5a},  we first  show that the solution $(u_\alpha, v_\alpha)$  given by \eqref{b60} is a local
 minimizer of the functional $I_{\alpha\beta}$.

 \begin{lemma}\label{lem8}
  Let $(u_\alpha, v_\alpha)$ be  defined by \eqref{b60}. Then $(u_\alpha, v_\alpha)$   is a local minimizer of the functional $I_{\alpha\beta}$ in $\mathcal{A}$.
 \end{lemma}

 {\bf  Proof.} \quad   It is easy to check that for any $(u', v')\in
 \mathcal{A}$,
   \berr
    &&\partial_{c_1}I_{\alpha\beta}(u'+c_1(u', v'), v'+c_2(u', v')) \nm\\
    &&=2(\alpha+\beta)\left[\re^{2c_1}\ito\re^{2u_0+2u'}\ud x-\re^{c_1}Q_1(u', v', \re^{c_2})+\frac{2\pi n}{\alpha\beta}\right]=0,\\
    &&\partial_{c_2}I_{\alpha\beta}(u'+c_1(u', v'), v'+c_2(u', v')) \nm\\
    &&=2(\alpha+\beta)\left[\re^{2c_2}\ito\re^{2u_0+2u'}\ud x-\re^{c_2}Q_2(u', v', \re^{c_1})+\frac{2\gamma\pi n}{\alpha\beta}\right]=0,
   \eerr
and
 \berr
    &&\partial^2_{c^2_1}I_{\alpha\beta}(u'+c_1(u', v'), v'+c_2(u', v')) \nm\\
    &&=2(\alpha+\beta)\left[2\re^{2c_1}\ito\re^{2u_0+2u'}\ud x-\re^{c_1}Q_1(u', v', \re^{c_2})\right],\\
    &&\partial^2_{c^2_2}I_{\alpha\beta}(u'+c_1(u', v'), v'+c_2(u', v')) \nm\\
    &&=2(\alpha+\beta)\left[\re^{2c_2}\ito\re^{2u_0+2u'}\ud x-\re^{c_2}Q_2(u', v', \re^{c_1})\right],\\
    &&\partial^2_{c_1c_2}I_{\alpha\beta}(u'+c_1(u', v'), v'+c_2(u', v')) \nm\\
    &&=-2(\alpha+\beta)\gamma\re^{c_1}\re^{c_2}\ito\re^{u_0+u'+v'}\ud x.
\eerr

 Then, in view of  \eqref{b20} and \eqref{b21}, we obtain
 \berr
  &&\partial^2_{c^2_1}I_{\alpha\beta}(u_\alpha, v_\alpha) \nm\\
    &&=
    2(\alpha+\beta)\left\{\left[(1-\gamma)\ito\re^{u_0+u_\alpha}\ud x+\gamma\ito\re^{u_0+u_\alpha+v_\alpha}\ud x\right]^2-\frac{8\pi n}{\alpha\beta}\ito\re^{2u_0+2u_\alpha}\ud x\right\}^{\frac12},\\
     &&\partial^2_{c^2_2}I_{\alpha\beta}(u_\alpha, v_\alpha) \nm\\
    &&=
    2(\alpha+\beta)\left\{\left[(1-\gamma)\ito\re^{v_\alpha}\ud x+\gamma\ito\re^{u_0+u_\alpha+v_\alpha}\ud x\right]^2-\frac{8\gamma\pi n}{\alpha\beta}\ito\re^{2v_\alpha}\ud x\right\}^{\frac12}.
 \eerr

 Since $(u_\alpha', v_\alpha')$ lies in the interior of $\mathcal{A}$, we  obtain
 \berr
   \partial^2_{c^2_1}I_{\alpha\beta}(u_\alpha, v_\alpha)>2(\alpha+\beta)\gamma\ito\re^{u_0+u_\alpha+v_\alpha}\ud x,\\
    \partial^2_{c^2_2}I_{\alpha\beta}(u_\alpha, v_\alpha)>2(\alpha+\beta)\gamma\ito\re^{u_0+u_\alpha+v_\alpha}\ud x.
  \eerr
 Hence,  at the point $(u_\alpha, v_\alpha)$ the Hessian matrix  of $I_{\alpha\beta}(u'+c_1, v'+c_2)$ with respect to $(c_1, c_2)$
 is strictly positive definite. Let $u=u'+c_1, v=v'+c_2$. By  continuity, there exists  $\delta>0$ such that for
  \[ \|u-u_\alpha\|+\|v-v_\alpha\|\le \delta, \]
 we have $(u', v')$ lies in the interior of $\mathcal{A}$ and
  \berr
   I_{\alpha\beta}(u, v)\ge I_{\alpha\beta}(u'+c_1(u', v'), v'+c_2(u',  v'))=J_{\alpha\beta}(u', v').
  \eerr

   Hence we have
 \berr
   I_{\alpha\beta}(u, v)\ge J_{\alpha\beta}(u', v')
   \ge \inf\limits_{(u', v')\in\mathcal{A}}J_{\alpha\beta}(u', v')=I_{\alpha\beta}(u_\alpha, v_\alpha),
  \eerr
 which says that $(u_\alpha, v_\alpha)$ is a local minimizer of $ I_{\alpha\beta}(u, v)$.

\subsection{The second solution}

 To find a secondary solution which is actually a mountain-pass critical point,  we show that the functional
 $I_{\alpha\beta}$ satisfies the  Palais--Smale condition.
  \begin{lemma}\label{lem9}
  Let $\{(u_m, v_m)\}$ be a sequence in $W^{1, 2}(\Omega)\times W^{1,  2}(\Omega)$ satisfying
   \ber
    I_{\alpha\beta}(u_m, v_m)\to \nu \quad   \text{as}\quad m\to +\infty,\label{b65}\\
    \|DI_{\alpha\beta}(u_m, v_m)\|_*\to 0 \quad   \text{as}\quad m\to +\infty,\label{b66}
   \eer
 where $\nu$ is a constant, $\|\cdot\|_*$ denotes the norm of the dual space of $W^{1,2}(\Omega)\times W^{1,2}(\Omega)$. Then $\{(u_m, v_m)\}$ admits a convergent subsequence in  $W^{1, 2}(\Omega)\times W^{1,  2}(\Omega)$ .

  \end{lemma}

 {\bf Proof.}  Let $u_m=c_{1, m}+u_m', v_m=c_{2, m}+v_m',$ where
  \[\ito u_m'\ud x=0,\quad  \ito v_m'\ud x=0,\quad c_{1, m}=\frac{1}{|\Omega|}\ito u_m\ud x, \quad c_{2, m}=\frac{1}{|\Omega|}\ito v_m\ud x.\]
  A simple computation gives
   \ber
  &&(DI_{\alpha\beta}(u_m, v_m))(\varphi, \psi)\nm\\
  &&= \frac12\left(\frac1\alpha+\frac1\beta\right)\ito(\nabla u_m\cdot\nabla\varphi+\nabla v_m\cdot\nabla\psi)\ud x
  +\frac12\left(\frac1\alpha-\frac1\beta\right)\ito(\nabla v_m\cdot\nabla\varphi+\nabla u_m\cdot\nabla\psi)\ud x \nm\\
  &&\quad +2\alpha\ito(\re^{u_0+u_m}+\re^{v_m}-2)\re^{u_0+u_m}\varphi\ud x+2\alpha\ito(\re^{u_0+u_m}+\re^{v_m}-2)\re^{v_m}\psi\ud x\nm\\
  &&\quad +2\beta\ito(\re^{u_0+u_m}-\re^{v_m})\re^{u_0+u_m}\varphi\ud x-2\beta\ito(\re^{u_0+u_m}-\re^{v_m})\re^{v_m}\psi\ud x \nm\\
  &&\quad +\frac{4\pi  n}{|\Omega|}\left(\frac1\alpha+\frac1\beta\right)\ito\varphi\ud x+\frac{4\pi  n}{|\Omega|}\left(\frac1\alpha-\frac1\beta\right)\ito\psi\ud x.\label{b67}
   \eer
 Taking $(\varphi, \psi)=(1, 0)$ and  $(\varphi, \psi)=(0, 1)$ in
 \eqref{b67},  from \eqref{b66}  we obtain as $m\to +\infty$,
   \ber
   (\alpha+\beta)\ito\re^{2u_0+2u_m}\ud x-2\alpha\ito\re^{u_0+u_m}\ud x-(\beta-\alpha)\ito\re^{u_0+u_m+v_m}\ud x+ 2\pi n\left(\frac1\alpha+\frac1\beta\right)\to 0,\label{b68}\\
   (\alpha+\beta)\ito\re^{2v_m}\ud x-2\alpha\ito\re^{v_m}\ud x-(\beta-\alpha)\ito\re^{u_0+u_m+v_m}\ud x+ 2\pi n\left(\frac1\alpha-\frac1\beta\right)\to  0.\quad \label{b69}
   \eer
 Combining \eqref{b68} and \eqref{b69} gives
  \ber
  (\alpha+\beta)\ito\left(\re^{2u_0+2u_m}+\re^{2v_m}\right)\ud x-2(\beta-\alpha)\ito\re^{u_0+u_m+v_m}\ud x\nm\\
  -2\alpha\ito(\re^{u_0+u_m}+\re^{v_m})\ud x+\frac{4\pi n}{\alpha}\to 0 \label{b70}
  \eer
   as $m\to +\infty$.

Noting that
  \berr
  &&\alpha\left(\re^{u_0+u_m}+\re^{v_m}-2\right)^2+\beta\left(\re^{u_0+u_m}-\re^{v_m}\right)^2\nm\\
  &&=(\alpha+\beta)\left(\re^{2u_0+2u_m}+\re^{2v_m}\right)-2(\beta-\alpha)\re^{u_0+u_m+v_m}-4\alpha(\re^{u_0+u_m}+\re^{v_m})+4\alpha
  \eerr
and
 \berr
 \alpha\left(\re^{u_0+u_m}+\re^{v_m}-2\right)^2+\beta\left(\re^{u_0+u_m}-\re^{v_m}\right)^2
 &\ge&\alpha\left(\left[\re^{u_0+u_m}+\re^{v_m}-2\right]^2+\left[\re^{u_0+u_m}-\re^{v_m}\right]^2\right)\nm\\
  &=&2\alpha\left(\left[\re^{u_0+u_m}-1\right]^2+\left[\re^{v_m}-1\right]^2\right),
 \eerr
 by \eqref{b70} we have
 \ber
  2\alpha\ito\left[\left(\re^{u_0+u_m}-1\right)^2+\left(\re^{v_m}-1\right)^2\right]\ud x
   +2\alpha\ito\left(\re^{u_0+u_m}+\re^{v_m}\right)\ud x+\frac{4\pi n}{\alpha}<4|\Omega|\alpha+o(1) \label{b71}
 \eer
 as $m\to +\infty$.

 Then it follows from \eqref{b71} that
  \ber
   \ito\left(\re^{u_0+u_m}+\re^{v_m}-2\right)^2\ud x+\ito\left(\re^{u_0+u_m}-\re^{v_m}\right)^2\ud x \le
  4|\Omega|+o(1),\label{b72}\\
  \ito\left(\re^{u_0+u_m}-1\right)^2\ud x+\ito\left(\re^{v_m}-1\right)^2\ud x
   \le 2|\Omega|+o(1),\label{b73}\\
    \ito \re^{u_0+u_m}\ud x\le2|\Omega|+o(1), \quad  \ito\re^{v_m}\ud x\le 2|\Omega|+o(1),\label{b74}
  \eer
 as  $m\to +\infty$.

 Using \eqref{b74} and the  Jensen inequality we obtain
   \ber
   \re^{c_{1,m}}\le 2+o(1), \quad \re^{c_{2,m}}\le 2+o(1),\label{b75}
   \eer
 as  $m\to +\infty$, which says $c_{1, m}, c_{2, m}$ are bounded from above.  From \eqref{b73}--\eqref{b74} we see that
 \ber
 \ito\re^{2u_0+2u_m}\ud x\le 6|\Omega|+o(1), \quad \ito\re^{2v_m}\ud x\le 6|\Omega|+o(1),\label{b76}
 \eer
 as  $m\to +\infty$.

Denote  $(u_m'+v_m')^+\equiv\max\{u_m'+v_m', 0\}$.  Setting   $\varphi=\psi=(u_m'+v_m')^+$ in \eqref{b67}, we have
  \ber
  &&(DI_{\alpha\beta}(u_m, v_m))\left[(u_m'+v_m')^+, (u_m'+v_m')^+\right]\nm\\
  &&=\frac1\alpha\|\nabla (u_m'+v_m')^+\|_2^2+2(\alpha+\beta)\ito\left(\re^{u_0+u_m}-\re^{v_m}\right)^2(u_m'+v_m')^+\ud x\nm\\
  &&\quad+8\alpha\ito\re^{u_0+u_m+v_m}(u_m'+v_m')^+\ud x-4\alpha\ito\left(\re^{u_0+u_m}+\re^{v_m}\right)(u_m'+v_m')^+\ud x\nm\\
  &&\quad +\frac{8\pi n}{\alpha|\Omega|}\ito(u_m'+v_m')^+\ud x\nm\\
  &&\ge8\alpha\ito\re^{u_0+u_m+v_m}(u_m'+v_m')^+\ud x-4\alpha\ito\left(\re^{u_0+u_m}+\re^{v_m}\right)(u_m'+v_m')^+\ud x\nm\\
  &&\ge8\alpha\ito\re^{u_0+u_m+v_m}(u_m'+v_m')^+\ud x\nm\\
  &&\quad -4\alpha\left(\left[\ito\re^{2u_0+2u_m}\ud x\right]^{\frac12}+\left[\ito\re^{2v_m}\ud x\right]^{\frac12}\right)\|(u_m'+v_m')^+\|_2.\label{b77}
 \eer

Then it follows from \eqref{b66}, \eqref{b76} and \eqref{b77} that
 \ber
  \ito\re^{u_0+u_m+v_m}(u_m'+v_m')^+\ud x&\le& C
  \left(\|(u_m'+v_m')^+\|_2+\vep_n\|(u_m'+v_m')^+\|\right)\nm\\
  &\le& C\left(\|\nabla u_m'\|_2+\|\nabla v_m'\|_2\right),\label{b78}
  \eer
where $C$ is a suitable positive constant  independent of $m$ and
$\vep_m\to 0$ as $m\to \infty$.

Now let $(\varphi, \psi)=(u_m', v_m')$ in \eqref{b67}, we see that
 \ber
 &&(DI_{\alpha\beta}(u_m, v_m))(u_m', v_m')\nm\\
 &&\ge \frac1\beta\left(\|\nabla u_m'\|_2^2+\|\nabla v_m'\|_2^2\right)+2(\alpha+\beta)\ito\re^{2u_0+2u_m}u_m'\ud x\nm\\
 && \quad-4\alpha\ito\re^{u_0+u_m}u_m'\ud x+2(\alpha+\beta)\ito\re^{2v_m}v_m'\ud x-4\alpha\ito\re^{v_m}v_m'\ud x\nm\\
 && \quad-2(\beta-\alpha)\ito\re^{u_0+u_m+v_m}(u_m'+v_m')^+\ud x \nm\\
 &&\ge \frac1\beta\left(\|\nabla u_m'\|_2^2+\|\nabla v_m'\|_2^2\right)+2(\alpha+\beta)\left[\ito\re^{2u_0}\re^{2c_{1, m}}(\re^{2u_m'}-1)u_m'\ud x+\ito\re^{2u_0}\re^{2c_{1, m}}u_m'\ud x\right]\nm\\
 && \quad +2(\alpha+\beta)\left[\ito\re^{2c_{2, m}}(\re^{2v_m'}-1)v_m'\ud x+\ito\re^{2c_{2, m}}v_m'\ud x\right]\nm\\
 && \quad-4\alpha\left[\left(\ito\re^{2u_0+2u_m}\ud x\right)^{\frac12}\|\nabla u_m'\|_2+\left(\ito\re^{2v_m}\ud x\right)^{\frac12}\|\nabla v_m'\|_2\right]-C\left(\|\nabla u_m'\|_2+\|\nabla v_m'\|_2\right)\nm\\
 &&\ge \frac1\beta\left(\|\nabla u_m'\|_2^2+\|\nabla v_m'\|_2^2\right)-C\left(\|\nabla u_m'\|_2+\|\nabla v_m'\|_2\right),\label{b79}
 \eer
 where we have used  \eqref{b75}--\eqref{b78} and the  H\"{o}lder inequality, $C$ is a positive constant independent of $m$.

Then using  \eqref{b66}--\eqref{b79}, we conclude  that there exits a  positive constant $C$ independent of $m$ such that
  \be
  \|\nabla u_m'\|_2+\|\nabla v_m'\|_2\le C.\label{b80}
  \ee

On the other hand, from \eqref{b65} we see that $I_{\alpha\beta}(u_m, v_m)$ is bounded from below. Therefore we have
  \ber
  && 4\pi n\left[\left(\frac1\alpha+\frac1\beta\right)c_{1, m}+\left(\frac1\alpha-\frac1\beta\right)c_{2, m}\right]\nm\\
  &&\ge -C-\alpha \ito\left(\re^{u_0+u_m}+\re^{v_m}-2\right)^2\ud x-\beta\ito\left(\re^{u_0+u_m}-\re^{v_m}\right)^2\ud x\nm\\
  &&\quad -\frac{1}{2\alpha}\left(\|\nabla u_m'\|_2^2+\|\nabla v_m'\|_2^2\right)\nm\\
  &&\ge -C, \label{b81}
  \eer
where we have used \eqref{b72} and \eqref{b80}, $C$ is a positive constant independent of $m$. Then  we infer from  \eqref{b75} and \eqref{b81} that $c_{1, m}, c_{2, m}$
are bounded from below uniformly with respect to $m$. As a result  $c_{1, m}, c_{2,  m}$
are bounded for all $m$.  Therefore we obtain  that $u_m=u_m'+c_{1, m}, v_m=v_m'+c_{2,  m}$ are bounded sequence in $W^{1, 2}(\Omega)$.
Hence there exists a subsequence of $u_m$ and $v_m$, still denoted by $u_m$ and $v_m$, such that
    \[
    u_m\to \tilde{u}, \quad v_m\to \tilde{v}
    \]
 weakly in $W^{1, 2}(\Omega)$, strongly in $L^p(\Omega)$ for any $p\ge1$ and pointwise a.e. in $\Omega$   as $m\to \infty$,
 where $\tilde{u}, \tilde{v}\in W^{1, 2}(\Omega)$.  Moreover, we  see that
   \berr
   \re^{u_m}\to \re^{\tilde{u}}, \quad   \re^{v_m}\to \re^{\tilde{v}} \quad \text{strongly in }\quad L^p(\Omega), \,\forall\, p\ge1
   \eerr
and
 \berr
 c_{1, m}=\frac{1}{|\Omega|}\ito u_m\ud x\to \tilde{c}_1=\frac{1}{|\Omega|}\ito \tilde{u}\ud x, \quad c_{2, m}=\frac{1}{|\Omega|}\ito v_m\ud x\to \tilde{c}_2=\frac{1}{|\Omega|}\ito \tilde{v}\ud x
 \eerr
 as $m\to \infty$.

 Hence, for any  $(\varphi, \psi)\in W^{1, 2}(\Omega)\times W^{1,  2}(\Omega)$, as $m\to \infty$  we have
  \be
  I'_{\alpha\beta}(u_m, v_m)(\varphi, \psi)\to I'_{\alpha\beta}(\tilde{u}, \tilde{v})(\varphi, \psi)=0, \label{b82}
  \ee
which says that $(\tilde{u}, \tilde{v})$ is a critical point for the
functional $I_{\alpha\beta}$.

 In the sequel we show the strong convergence of $(u_m, v_m)$.
  Let $(\varphi, \psi)=(u_m-\tilde{u}, v_m-\tilde{v})$. Then, from
  \eqref{b66} and \eqref{b82},  we see that
   \ber
   [DI_{\alpha\beta}(u_m, v_m)- DI_{\alpha\beta}(\tilde{u}, \tilde{v})][u_m-\tilde{u}, v_m-\tilde{v}]\le\vep_n\left(\|u_m-\tilde{u}\|+\|v_m-\tilde{v}\|\right)=o(1),\label{b83}
   \eer
    as $m\to \infty$.

 Then from \eqref{b83} we conclude that
  \berr
  && \frac1\beta\left(\|\nabla (u_m-\tilde{u})\|_2^2+\|\nabla (v_m-\tilde{v})\|_2^2\right)\nm \\
  &&\le -2(\alpha+\beta)\left(\ito\re^{2u_0}[\re^{2u_m}-\re^{2\tilde{u}}][u_m-\tilde{u}]\ud x+\ito[\re^{2v_m}-\re^{2\tilde{v}}][v_m-\tilde{v}]\ud x\right)\nm \\
  &&\quad +2(\alpha-\beta)\ito\re^{u_0}(\re^{u_m+v_m}-\re^{\tilde{u}+\tilde{v}})([u_m-\tilde{u}]+[v_m-\tilde{v}])\ud x\nm\\
  &&\quad +4\alpha\ito \left(\re^{u_0}[\re^{u_m}-\re^{\tilde{u}}][u_m-\tilde{u}]+[\re^{v_m}-\re^{\tilde{v}}][v_m-\tilde{v}]\right)\ud x+o(1)\nm\\
  &&=o(1),
  \eerr
 as $m\to \infty$. Therefore we see that $(u_m, v_m)$ converges strongly to $(\tilde{u}, \tilde{v})$ in $W^{1,2}(\Omega)\times W^{1,2}(\Omega)$ as $m\to \infty$. Then the lemma follows.

 We have shown that, for $\alpha$ large enough, there exists a solution $(u_\alpha, v_\alpha)$ of \eqref{b4a}--\eqref{b5a}  satisfying the behavior
 \[\re^{u_0+u_\alpha}\to1, \quad \re^{v_\alpha}\to1,\]
 as $\alpha\to+\infty$ and this solution is a local minimum for $I_{\alpha\beta}$ by Lemma \ref{lem8}. Then there exists a positive constant $\delta_0>0$ such that
  \[ I_{\alpha\beta}(u, v)\ge I_{\alpha\beta}(u_\alpha, v_\alpha)\quad  \text{if}\quad \|u-u_\alpha\|+\|v-v_\alpha\|\le \delta_0.\]

To find a secondary solution, we note that the functional admits a mountain-pass structure. Indeed,   noting that $u_0+u_\alpha<0, v_\alpha<0$,   we have for all $c<0$,
 \ber
  I_{\alpha\beta}(u_\alpha+c, v_\alpha)- I_{\alpha\beta}(u_\alpha, v_\alpha)\le (4\alpha+\beta)|\Omega|+4\pi n\left(\frac1\alpha+\frac1\beta\right)c.
  \label{b84}
 \eer

We consider two situations.

(1) If $(u_\alpha, v_\alpha)$ is not a strictly local minimum
for  the functional $I_{\alpha\beta}$, i.e., for all
$0<\delta<\delta_0$,
 \berr
  \inf\limits_{\|u-u_\alpha\|+\|v-v_\alpha\|= \delta}  I_{\alpha\beta}(u, v)
   =I_{\alpha\beta}(u_\alpha, v_\alpha).
 \eerr
Then by using Ekeland's lemma (see \cite{gho}), we get a local
minimum
 $(u_\alpha^\delta, v_\alpha^\delta)$ for $I_{\alpha\beta}$ such
 that
  \[\|u_\alpha^\delta-u_\alpha\|+\|v_\alpha^\delta-v_\alpha\|= \delta \quad \text{for any}\quad \delta\in(0, \delta_0).\]
 Of course, we can find a secondary solution.

 (2) If $(u_\alpha, v_\alpha)$ is
a strictly local minimum for  the functional $I_{\alpha\beta}$, then
there exits $\delta_1\in(0, \delta_0)$ such that
 \[\inf\limits_{\|u-u_\alpha\|+\|v-v_\alpha\|= \delta_1}  I_{\alpha\beta}(u, v)
  >I_{\alpha\beta}(u_\alpha, v_\alpha). \]

  Moreover, in view of   \eqref{b84} we can choose $|\tilde{c}|>\delta_0$ such that
 \berr
   I_{\alpha\beta}(u_\alpha+\tilde{c}, v_\alpha)< I_{\alpha\beta}(u_\alpha, v_\alpha)-1
 \eerr

 Now we introduce the paths
  \[
   \mathcal{P}=\left\{\Gamma(t)|\Gamma\in C\big([0, 1];  W^{1, 2}(\Omega)\times W^{1, 2}(\Omega)\big):
   \,\, \Gamma(0)=(u_\alpha, v_\alpha), \, \, \Gamma(1)=(u_\alpha+\tilde{c}, v_\alpha)\right\}
  \]

  and define
  \[
  m_0=\inf\limits_{\Gamma\in\mathcal{P}}\sup\limits_{t\in[0,1]}\big\{I_{\alpha\beta}(\Gamma(t))\big\}.
  \]
Hence we have
 \be
 m_0>I_{\alpha\beta}(u_\alpha, v_\alpha).\label{b85}
 \ee
 Then by Lemma \ref{lem9} we see that the functional  $I_{\alpha\beta}$ satisfies all the hypotheses of the mountain-pass theorem of Ambrosetti--Rabinowitz \cite{amra}.
Therefore  we conclude that $m_0$ is a critical value of the functional $I_{\alpha\beta}$. Noting \eqref{b85}, we have an
additional solution of the equations \eqref{b4a} and \eqref{b5a}.

Then the proof of Theorem \ref{thb1} is complete.

\section{Summary and comments}

In this paper, we developed an  existence theory for the multiple vortex solutions of the non-Abelian Chern--Simons--Higgs equations in the Gudnason model \cite{gud1,gud2}
of ${\cal N}=2$ supersymmetric field theory where gauge field dynamics is governed by two Chern--Simons terms. In the full-plane situation, we proved the existence of solutions
for the general system of nonlinear partial differential equations involving an arbitrary number of unknowns and established the exponential decay estimates for the solutions which
realize a family of quantized integrals expressed in terms of vortex numbers. These solutions give rise to multiple vortex field configurations which approach  at spatial infinity the vacuum state with
completely broken symmetry and known as topological solutions. There are no restrictions to the values of various coupling parameters and vortex numbers
or vortex distributions for the existence results to hold. In the doubly periodic domain setting, we confined our study to the case when the vortex equations contain only two unknowns.
We derived a necessary condition for the existence of a multiple vortex solution and obtained some sufficient conditions for the existence of two distinct solutions for a given prescribed
distribution of vortices. We also established the asymptotic behavior of the solutions as a coupling parameter goes to infinity.

The problems which remain untouched and are of considerable future interest are the existence of nontopological solutions realizing the asymptotic vacuum state with unbroken symmetry
and the existence of solutions of the multiple vortex equations involving more than two unknowns. Our methods developed so far are still not sufficiently effective in dealing with these
rather challenging problems and new ideas and techniques are called upon in order to make further progress in this area.
\small{

}
\end{document}